\documentclass[12pt]{article}

\usepackage{amsmath}
\usepackage{amsfonts}
\usepackage{latexsym}
\usepackage{graphicx}

\newtheorem{proposition}{Proposition}[section]

\newtheorem{remark}{Remark}[section]

\usepackage{xcolor}

\date{ }

\begin{document}

\title{Topological optimization and minimal compliance in linear elasticity}

\author{Cornel Marius Murea$^1$, Dan Tiba$^2$\\
{\normalsize $^1$ D\'epartement de Math\'ematiques, IRIMAS,}\\
{\normalsize Universit\'e de Haute Alsace, France,}\\
{\normalsize cornel.murea@uha.fr}\\
{\normalsize $^2$ Institute of Mathematics (Romanian Academy) and}\\ 
{\normalsize Academy of Romanian Scientists, Bucharest, Romania,}\\ 
{\normalsize dan.tiba@imar.ro}
}

\maketitle

\begin{abstract}
We investigate a fixed domain approach in shape optimization, using a regularization
of the Heaviside function both in the cost functional and in the state system.
We consider the compliance minimization problem in linear elasticity, a well known
application in this area of research. The optimal design problem is approached by
an optimal control problem defined in a prescribed domain including all the admissible unknown domains.
This approximating optimization problem has good differentiability properties and
a gradient algorithm can be applied. Moreover, the paper also includes several numerical
experiments that demonstrate the descent of the obtained cost values and show the
topological and the boundary variations of the computed domains.
{\color{black}The proposed approximation technique is new and can be applied to state
  systems given by various boundary value problems. }

\textbf{Keywords: topological optimization; minimal compliance}
\end{abstract}

\section{Minimal compliance in linear elasticity}

Let $\Omega\subset \mathbb{R}^2$ be a connected, bounded, Lipschitz domain with boundary
{\color{black}$\partial \Omega
  =\overline{\Gamma}_D \cup \overline{\Gamma}_N \cup \overline{\Gamma}$,
where $\Gamma _D$, $\Gamma _N$ and $\Gamma$ are relatively open subsets, mutually disjoint,
}
such that $meas(\Gamma_D)>0$, $meas(\Gamma_N)>0$, $meas(\Gamma)>0$.

The notation $\mathbf{v}\cdot \mathbf{w}$ means the scalar product of two vectors
$\mathbf{v}, \mathbf{w}\in \mathbb{R}^2$ and $A:B =\sum_{i,j=1}^2 a_{ij}b_{ij}$ if
$A=(a_{ij})_{1\leq i,j\leq 2}$, $B=(b_{ij})_{1\leq i,j\leq 2}$.
We will use
$$
\nabla \cdot \mathbf{v}
=\frac{\partial v_1}{\partial x_1}+\frac{\partial v_2}{\partial x_2},\quad
\nabla \cdot A
=\left(
\begin{array}{c}
  \frac{\partial a_{11}}{\partial x_1}+\frac{\partial a_{12}}{\partial x_2}\\
  \frac{\partial a_{21}}{\partial x_1}+\frac{\partial a_{22}}{\partial x_2}
\end{array}
\right)
$$
to denote the divergence operator of a vector valued function $\mathbf{v}\in \mathbb{R}^2$ or
 $A=(a_{ij})_{1\leq i,j\leq 2}\in \mathbb{R}^{2\times 2}$.

We denote by $\mathbf{y}:\overline{\Omega}\rightarrow \mathbb{R}^2$ the displacement of
the linear elastic body $\Omega\subset \mathbb{R}^2$.
The stress tensor in linear elasticity is given by
$$
\sigma \left(\mathbf{y}\right) = 
\lambda^S(\nabla \cdot \mathbf{y})\mathbf{I}+2\mu^S \mathbf{e}(\mathbf{y})
$$
where $\lambda^S,\mu^S>0$ are the Lam\'e coefficients,
{\color{black}independent of the space variable, }
$\mathbf{I}$ is the unity matrix
and $\mathbf{e}(\mathbf{y})=\frac{1}{2}\left(\nabla\mathbf{y} + (\nabla\mathbf{y})^T\right)$.

For given volume load $\mathbf{f}:\Omega\rightarrow \mathbb{R}^2$ and
surface load $\mathbf{h}:\Gamma_N\rightarrow \mathbb{R}^2$,
we consider the linear elasticity equations:
find $\mathbf{y}:\overline{\Omega}\rightarrow \mathbb{R}^2$ such that
\begin{eqnarray}
  -\nabla\cdot \sigma \left(\mathbf{y}\right) & = & \mathbf{f},\hbox{ in }\Omega
  \label{elast1}\\
  \mathbf{y} & = & 0,\hbox{ on }\Gamma_D\label{elast2}\\
  \sigma \left(\mathbf{y}\right)\mathbf{n}& = & \mathbf{h},\hbox{ on }\Gamma_N\label{elast3}\\
  \sigma \left(\mathbf{y}\right)\mathbf{n}& = & 0,\hbox{ on }\Gamma\label{elast4}
\end{eqnarray}
where $\mathbf{n}$ is the unit outer normal vector along the boundary.

The weak formulation is: find $\mathbf{y}\in V$ such that
\begin{equation}
  \int_\Omega  \sigma \left( \mathbf{y}\right) : \nabla \mathbf{v}\,d\mathbf{x}
   = 
  \int_\Omega \mathbf{f}\cdot \mathbf{v}\,d\mathbf{x}
+\int_{\Gamma_N} \mathbf{h}\cdot \mathbf{v}\,ds,\quad \forall \mathbf{v} \in V
\label{weakelast}
\end{equation}
where $\mathbf{f}\in \left(L^2(\Omega)\right)^2$,
$\mathbf{h}\in \left(L^2(\Gamma_N)\right)^2$,
$$
V=\{\mathbf{v} \in \left(H^1(\Omega)\right)^2;\ \mathbf{v}=0\hbox{ on }\Gamma_D \}.
$$
It is well known that this problem has a unique solution in $V$, see \cite{Ciarlet1988}.

Using that $\sigma \left(\mathbf{y}\right)$ is symmetric and the identity
$A:B=A^T:B^T$, we obtain that the left-hand side from (\ref{weakelast})
with $\mathbf{v} = \mathbf{y}$ satisfies:
\begin{eqnarray}
&&\int_\Omega  \sigma \left( \mathbf{y}\right) : \nabla \mathbf{y}\,d\mathbf{x}
=\int_\Omega  \frac{1}{2}\sigma \left( \mathbf{y}\right) : \nabla \mathbf{y}
+\frac{1}{2}\sigma \left( \mathbf{y}\right) : \nabla \mathbf{y}\,d\mathbf{x}\nonumber\\
&&
=\int_\Omega  \frac{1}{2}\sigma \left( \mathbf{y}\right) : \nabla \mathbf{y}
+\frac{1}{2}\left(\sigma \left( \mathbf{y}\right)\right)^T
: \left(\nabla \mathbf{y}\right)^T\,d\mathbf{x}\nonumber\\
&&
=\int_\Omega  \frac{1}{2}\sigma \left( \mathbf{y}\right) : \nabla \mathbf{y}
+\frac{1}{2}\sigma \left( \mathbf{y}\right)
: \left(\nabla \mathbf{y}\right)^T\,d\mathbf{x}
\nonumber\\
&&
=\int_\Omega  \sigma \left( \mathbf{y}\right) : \mathbf{e}(\mathbf{y})\,d\mathbf{x}
=\int_\Omega \lambda^S(\nabla \cdot \mathbf{y})\mathbf{I}: \mathbf{e}(\mathbf{y})
+2\mu^S\mathbf{e}(\mathbf{y}):\mathbf{e}(\mathbf{y})\,d\mathbf{x}\nonumber\\
&&
=\int_\Omega \lambda^S(\nabla \cdot \mathbf{y})^2
+2\mu^S\mathbf{e}(\mathbf{y}):\mathbf{e}(\mathbf{y})\,d\mathbf{x}.
\label{coercive}
\end{eqnarray}

\begin{figure}[ht]
\begin{center}
\includegraphics[width=6cm]{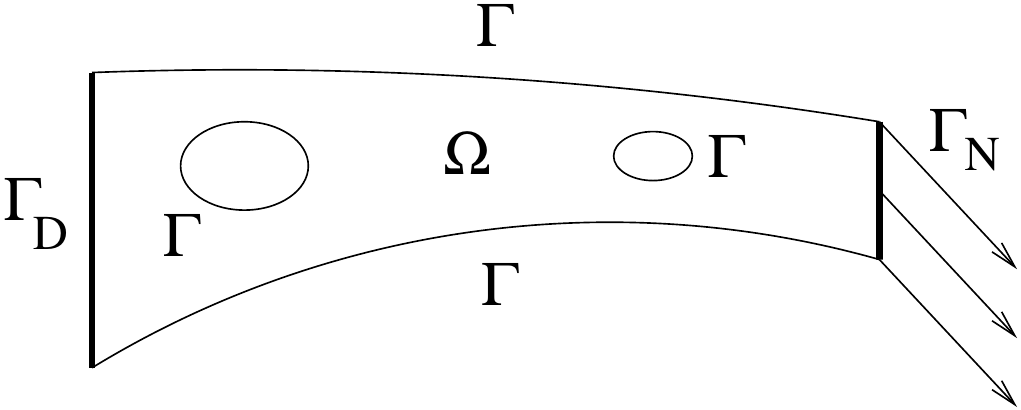}
\end{center}
\caption{The geometrical configuration.}
\label{fig:omega}
\end{figure}

A classical problem in structural design, see \cite{Bendsoe1995},
\cite{Bendsoe2003}, \cite{Allaire2007}, is to find
a domain $\Omega$ that minimizes the compliance 
(the work done by the load, 
expressed by the right-hand side in (\ref{weakelast}) with $\mathbf{v} = \mathbf{y}$) 
subject to
$\Gamma_N\subset \partial \Omega$,
$\Gamma_D \subset \partial \Omega$ and
the volume of $\Omega$ is prescribed, see Figure \ref{fig:omega}.
We suppose that $\Gamma_N$ and $\Gamma_D$ are fixed.
On $\Gamma_D$ the elastic body $\Omega$ is also fixed due to (\ref{elast2}),
while $\Gamma_N$ is fixed {\color{black} in the sense that this part of the boundary is
specified in advance for the family of all admissible domains $\Omega$
and may deform under the action of the traction $\mathbf{h}$
and the volume load $\mathbf{f}$ in (\ref{elast3}). }

Here, we examine the minimization of the compliance as well as of the volume of $\Omega$.
In practice, we penalize the volume of $\Omega$, and
the function to minimize is
\begin{equation}
\int_\Omega \mathbf{f}\cdot \mathbf{y}\,d\mathbf{x}
+\int_{\Gamma_N} \mathbf{h}\cdot \mathbf{y}\,ds
+\ell \int_\Omega 1\,d\mathbf{x}
\label{compliance1}
\end{equation}
where $\ell >0$ is a penalization coefficient.
The fixed domain method that we introduce here combines boundary variations with topology
optimization (the domain $\Omega$ is not necessarily simply connected and the number
of holes may change during the iterations). This is characterized by $\Gamma$, the
part of $\partial \Omega$ that is not fixed.
{\color{black}It is a new approach and another important property is that it may be applied to
  many boundary value problems as governing systems. }
For other fixed domain approaches, we quote \cite{Tiba1992}, \cite{Tiba2009}, \cite{Tiba2018},
\cite{Tiba2018a} and the survey \cite{Tiba2012} with its references.
For multi-layered composite materials, one can consult \cite{Delgado2017}.

\section{The shape optimization problem in fixed domain and its gradient}

\begin{figure}[ht]
\begin{center}
\includegraphics[width=6cm]{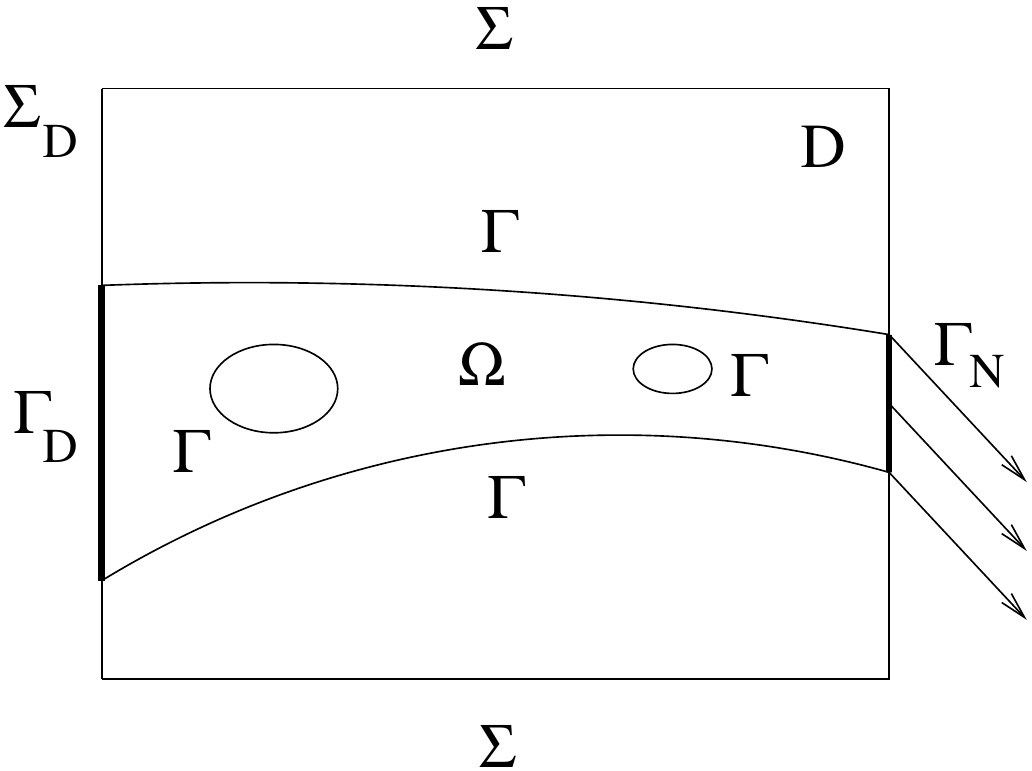}
\end{center}
\caption{The fixed domain $D$ including the unknown domain.}
\label{fig:D}
\end{figure}

We consider a simply connected, bounded, Lipschitz domain $D \subset \mathbb{R}^2$, including
the unknown domain $\Omega$, with
{\color{black}$\partial D=\overline{\Sigma}_D \cup \overline{\Gamma}_N \cup \overline{\Sigma}$,
where $\Sigma_D$, $\Gamma _N$ and $\Sigma$
are relatively open subsets, mutually disjoint,
}
 such that $\Gamma_D \subset \Sigma_D$, see Figure \ref{fig:D}.

Let $X(D)$ denote a cone of $\mathcal{C}(\overline{D})$.
Following \cite{Tiba2009}, \cite{Tiba2012}, with any $g\in X(D)$, that we
call a parametrization, we associate the open set
\begin{equation*}
\Omega_g = int \left\{\mathbf{x} \in D;\ g(\mathbf{x})\geq 0\right\}.
\end{equation*}

{\color{black}
We define the family of admissible domains as the connected components of all
$\Omega_g, \; g \in X(D)$ satisfying  $\Gamma_N\subset \partial \Omega_g$,
$\Gamma_D \subset \partial \Omega_g$. This family is very rich and some examples are commented
in Remark \ref{rem:3.1} below. }

We use the following
regularization of the Heaviside function
\begin{equation}\label{He}
H^\epsilon(r)=
\left\{
\begin{array}{ll}
1-\frac{1}{2}e^{-\frac{r}{\epsilon}}, & r \geq 0,\\
\frac{1}{2}e^{\frac{r}{\epsilon}}, & r < 0 ,
\end{array}
\right.
\end{equation}
{\color{black}where $\epsilon >0$ is a parameter. }
We have that $H^\epsilon(g)$ is a regularization of the
characteristic function of $\overline{\Omega}_g$. This procedure has already been introduced
in \cite{Tiba1992}. Another approximation that we shall also use is
$H_\epsilon(r) =1,\ r \geq 0$ and $H_\epsilon(r) = \epsilon,\ r<0$.
The advantage of (\ref{He}) is its differentiability.

For given $\mathbf{f}\in \left(L^2(D)\right)^2$, $\mathbf{h}\in \left(L^2(\Gamma_N)\right)^2$ and 
$\ell >0$, we introduce the control problem (with control $g$) that approximates the shape optimization problem
\begin{equation}\label{jeps}
\inf_{g\in X(D)} \{J(g)=
\int_D H^\epsilon(g) \mathbf{f}\cdot \mathbf{y}^\epsilon(g)\,d\mathbf{x}
+\int_{\Gamma_N} \mathbf{h}\cdot \mathbf{y}^\epsilon(g)\,ds
+\ell \int_D H^\epsilon(g) \,d\mathbf{x}\}
\end{equation}
where $\mathbf{y}^\epsilon(g)\in W$ is the solution of
\begin{eqnarray}
\int_D H^\epsilon(g) \sigma \left( \mathbf{y}^\epsilon(g)\right) : \nabla \mathbf{v}\,d\mathbf{x}
&=&\int_D H^\epsilon(g) \mathbf{f}\cdot \mathbf{v}\,d\mathbf{x}
+\int_{\Gamma_N} \mathbf{h}\cdot \mathbf{v}\,ds,
\label{stateeps}
\end{eqnarray}
for all $\mathbf{v} \in W$, where
$$
W=\{\mathbf{v} \in \left(H^1(D)\right)^2;\ \mathbf{v}=0\hbox{ on }\Gamma_D \}.
$$
{\color{black}
We point out that $H^\epsilon(g)>0$ in $D$ and $H^\epsilon(g) \geq 1/2$ in $\overline{\Omega}_g$ due
to (\ref{He}). Moreover, $H^\epsilon(g)$ is an approximation of the characteristic function of
$\overline{\Omega}_g$ and due to this, there is an approximation relation between the solutions
of (\ref{stateeps}) and (\ref{weakelast}), that we examine in the sequel. }
  
\medskip

{\color{black}
  Regularization methods have a long history in fixed domain methods for problems involving unknown
  domains. For instance, already \cite{kawarada} used a penalization/regularization method in free
  boundary problems. A survey on this subject is the paper \cite{Tiba2012} and in \cite{Tiba2018}, \cite{Tiba2013}
  such approaches are extended to the optimization of plates with holes and other problems. In general, Dirichlet
  boundary conditions are taken into account, while the approximation defined in
  (\ref{jeps}), (\ref{stateeps}) can be used for other boundary conditions as well. }

\begin{remark}\label{rem:3.1}
An admissible  domain $\Omega_g$ for $g\in \mathcal{C}(\overline{D})$
has to satisfy
$\Gamma_N \subset \partial \Omega_g$,
$\Gamma_D \subset \partial \Omega_g$, that can be obtained by imposing simple equality constraints on $g$ (see below). 
Moreover, we also assume $\Omega_g$ to be connected, and this constraint has to be added to
 the definition of $X(D)$. 
{\color{black}However in the regularized problem, this is not necessary since we work in $D$. }
Other constraints on the geometry $\Omega_g$ may be included in the 
definition below, if necessary. For instance $E  \subset \Omega_g$, where $E$ is some given domain 
(add below the inequality $g > 0$ in $E$), etc.

We set 
$$
X(D)=\left\{ g\in \mathcal{C}(\overline{D});\
g(\mathbf{x}) = 0,\ \mathbf{x}\in \Gamma_N\cup\Gamma_D 
\right\}
$$
which is a subspace in $\mathcal{C}(\overline{D})$.
If $g \in \mathcal{C}^1(D)$ and is noncritical on its null level set, then the condition from the beginning
of this remark is satisfied. In general, $X(D)$ is a cone.
\end{remark}

\begin{proposition}\label{prop:1}
The problem (\ref{stateeps}) has a unique solution $\mathbf{y}^\epsilon(g)\in W$
and
\begin{equation}\label{norme_yeps}
\left\| \mathbf{y}^\epsilon(g)\right\|_{1,D} \leq \frac{C}{c(\epsilon)}\left( 
\left\| \mathbf{f}\right\|_{0,D} + \left\| \mathbf{h}\right\|_{0,\Gamma_N}
\right)
\end{equation}
where $C>0$ is independent of $\epsilon$ and $c(\epsilon) >0$ is indicated below.
\end{proposition}

\noindent
\textbf{Proof.}
We set
$$
a(\mathbf{y},\mathbf{v}) 
=\int_D H^\epsilon(g) \sigma \left( \mathbf{y}\right) : \nabla \mathbf{v}\,d\mathbf{x}
$$
and we get from (\ref{coercive}):
\begin{eqnarray*}
a(\mathbf{v},\mathbf{v})
&=&
\int_D H^\epsilon(g) \sigma \left( \mathbf{v}\right) : \nabla \mathbf{v}\,d\mathbf{x}
\\
&=&
\int_D H^\epsilon(g) \left( \lambda^S(\nabla \cdot \mathbf{v})^2
+2\mu^S\mathbf{e}(\mathbf{v}):\mathbf{e}(\mathbf{v})\right) d\mathbf{x}
\\
&\geq&
c(\epsilon) \int_D \left( \lambda^S(\nabla \cdot \mathbf{v})^2
+2\mu^S\mathbf{e}(\mathbf{v}):\mathbf{e}(\mathbf{v})\right) d\mathbf{x}\\
&\geq&
c(\epsilon) \int_D 2\mu^S\mathbf{e}(\mathbf{v}):\mathbf{e}(\mathbf{v})\, d\mathbf{x}.
\end{eqnarray*}

We have used that $ H^\epsilon(g) \geq c(\epsilon) > 0$ in $D$, due to (\ref{He}). In fact $c(\epsilon) $ also depends on $g$, that is fixed here (and uniformly bounded).
From the Korn's inequality, see \cite{Ciarlet1988}, we obtain that
$a$ is $W$-elliptic, i.e. 
$a(\mathbf{v},\mathbf{v})\geq \frac{c(\epsilon)}{C}\left\| \mathbf{v}\right\|_{1,D}^2$, 
where $C>0$ is independent on $g$, $\epsilon$.
From the Lax-Milgram theorem, we get that the problem has a unique solution
and
$$
\left\| \mathbf{y}^\epsilon(g)\right\|_{1,D} \leq \frac{C}{c(\epsilon)}\left( 
\left\| H^\epsilon(g)\mathbf{f}\right\|_{0,D} + \left\| \mathbf{h}\right\|_{0,\Gamma_N}
\right)
\leq \frac{C}{c(\epsilon)}\left( 
\left\| \mathbf{f}\right\|_{0,D} + \left\| \mathbf{h}\right\|_{0,\Gamma_N}
\right).
$$
since $0 < H^\epsilon(g) \leq 1$.
\quad$\Box$

We indicate now some basic approximation results.

\begin{proposition}\label{propNEW:1}
When $\epsilon \rightarrow 0$, on a subsequence, we have
$\mathbf{y}^{\epsilon} |_{\Omega_g} \rightarrow \mathbf{y}$ weakly in $H^1 (\Omega_g)$.
Moreover, $\mathbf{y} \in V$ and satisfies (\ref{elast1}) - (\ref{elast4}) in the
{\color{black}distributional }
sense. {\color{black}This statement remains valid for the corresponding approximating solutions }
 when $H^\epsilon$ is replaced by $H_\epsilon$ in (\ref{stateeps}).
\end{proposition}

\noindent
\textbf{Proof.}
Let $\mathbf{f_1} \in L^2 (D)$ be the extension by $0$ of $\mathbf{f} \in L^2(\Omega_g)$.
We use the approximating formulation (\ref{stateeps}) with $\mathbf{f}$ replaced by
$\mathbf{f_1}$ and $\mathbf{v} = \mathbf{y}_{\epsilon} \in W$:
\begin{eqnarray}
\int_D H^\epsilon(g) \sigma \left( \mathbf{y}^\epsilon\right) :
\nabla \mathbf{y}^\epsilon\,d\mathbf{x}
&=&\int_D H^\epsilon(g) \mathbf{f_1}\cdot \mathbf{y}^{\epsilon}\,d\mathbf{x}
+\int_{\Gamma_N} \mathbf{h}\cdot \mathbf{y}^{\epsilon}\,ds.
\label{stateepsNEW}
\end{eqnarray}

By (\ref{coercive}) and (\ref{stateepsNEW}), we get
\begin{eqnarray}
&&\int_D H^\epsilon(g) [\lambda^S(\nabla \cdot \mathbf{y}^{\epsilon})^2
  +2\mu^S\mathbf{e}(\mathbf{y}^{\epsilon}):\mathbf{e}(\mathbf{y}^{\epsilon}) ]\,d\mathbf{x}
\nonumber\\
&&=\int_D H^\epsilon(g) \mathbf{f_1}\cdot \mathbf{y}^{\epsilon}\,d\mathbf{x}
+\int_{\Gamma_N} \mathbf{h}\cdot \mathbf{y}^{\epsilon}\,ds.
\label{estimateNEW}
\end{eqnarray}

In the left-hand side of (\ref{estimateNEW}), all the terms are positive,
due to (\ref{He}), and we infer the inequality
\begin{eqnarray}
\mu^S\int_{\Omega_g} \mathbf{e}(\mathbf{y}^{\epsilon}) :
\mathbf{e}(\mathbf{y}^{\epsilon}) \,d\mathbf{x} 
&\leq&\int_{\Omega_g}  \mathbf{f}\cdot \mathbf{y}^{\epsilon}\,d\mathbf{x}
+\int_{\Gamma_N} \mathbf{h}\cdot \mathbf{y}^{\epsilon}\,ds,
\label{estimateNEW2}
\end{eqnarray}
where we use that $1 \geq H^\epsilon(g) \geq 1/2$ in $\Omega_g$ and the definition of $\mathbf{f}_1$.
One can apply the Korn's inequality in (\ref{estimateNEW2}) and establish that
$\mathbf{y}^{\epsilon}|\Omega_g$ is bounded in $H^1 (\Omega_g)$. On a subsequence,
$\mathbf{y}^{\epsilon} |_{\Omega_g} \rightarrow \mathbf{y}$ weakly in $H^1 (\Omega_g)$.
Moreover, $ H^\epsilon(g) \rightarrow H(g)$ in $L^p(D)$, for any $p \geq 1$.

For any test function $\mathbf{v} \in \mathcal{D}(\Omega_g) \subset W$, we pass to the limit
in (\ref{stateeps})
and obtain
$$
\int_{\Omega_g} \sigma \left( \mathbf{y} \right) :
\nabla \mathbf{v} \,d\mathbf{x}
=\int_{\Omega_g}  \mathbf{f}\cdot \mathbf{v}\,d\mathbf{x}
+\int_{\Gamma_N} \mathbf{h}\cdot \mathbf{v}\,ds,
\ \forall \mathbf{v} \in \mathcal{D}(\Omega_g)
$$
and the proof of the first part is finished.

For the last statement of Proposition \ref{propNEW:1}, we denote by $\mathbf{y}_{\epsilon}$
the corresponding solution
and we assume that $\mathbf{f} \in L^2(D)$ is given (we don't use the above extension by $0$).
We use the definition of $H_\epsilon$ and decompose (\ref{stateeps}) as follows
\begin{eqnarray}
&&  \epsilon \int_D \sigma \left( \mathbf{y}_\epsilon\right) : \nabla \mathbf{v}\,d\mathbf{x}
+ (1-\epsilon) \int_{\Omega_g} \sigma \left( \mathbf{y}_\epsilon\right) :
\nabla \mathbf{v}\,d\mathbf{x}
\nonumber\\
&&=\epsilon \int_D \mathbf{f}\cdot \mathbf{v}\,d\mathbf{x} + (1-\epsilon) \int_{\Omega_g} \mathbf{f}\cdot \mathbf{v}\,d\mathbf{x}
+\int_{\Gamma_N} \mathbf{h}\cdot \mathbf{v}\,ds,\ \forall \mathbf{v} \in W.
\label{stateeps2}
\end{eqnarray}

Putting $\mathbf{v} = \mathbf{y}_{\epsilon}$ in (\ref{stateeps2}), it is easy to see,
by using the identity (\ref{coercive}) and Korn's inequality in $D$,
that the difference between the two terms with coefficient $\epsilon$ from (\ref{stateeps2})
is bounded from
below by some constant. The other terms are handled as above, using Korn's inequality in
$\Omega_g$ and we get that $\mathbf{y}_\epsilon$ is bounded in $H^1(\Omega_g)$ if
$\epsilon \leq 1/2$. The passage to the limit is as above.
\quad$\Box$

\begin{remark}\label{rem:3.new}
Concerning shape optimization problems and their approximation, it is advantageous to
use $H^\epsilon$ due to the differentiability properties that will be discussed below.
 When $g$ is the unknown control and may change, the  procedure to
extend $\mathbf{f}$ by $0$ outside the $supp(\mathbf{f})$ (from the first part of the proof) involves
the hypothesis that $g \geq 0$ in $supp(\mathbf{f})$, for any admissible $g$.  The numerical
examples from the last section confirm that our method allows topological and boundary
variations and ensures a good descent of the cost.
\end{remark}

\begin{proposition}\label{prop:2}
For any $g$, $w$ in $X(D)$, the mapping $g\rightarrow \mathbf{y}^\epsilon(g)\in W$ is G\^ateaux
differentiable at $g$ {\color{black}and the directional derivative in the direction $w$, }
denoted by $\mathbf{z}\in W$, is the unique solution of the problem
\begin{eqnarray}
&&\int_D H^\epsilon(g) \sigma \left( \mathbf{z}\right) : \nabla \mathbf{v}\,d\mathbf{x}=
\nonumber\\
&&
- \int_D (H^\epsilon)^\prime(g) w\, \sigma \left( \mathbf{y}^\epsilon(g)\right)
: \nabla \mathbf{v}\,d\mathbf{x}
+\int_D  (H^\epsilon)^\prime(g) w\, \mathbf{f}\cdot \mathbf{v}\,d\mathbf{x},\ \forall \mathbf{v} \in W.
\label{limit}
\end{eqnarray}
\end{proposition}

\noindent
\textbf{Proof.}
Let $g$, $w$ be fixed in $X(D)$ and $\lambda \neq 0$, small. 
We notice first that $X(D)$ is stable to small perturbations, i.e. $g + \lambda w \in X(D)$ if $|\lambda|$ small. 

We write the equation (\ref{stateeps}) for $g+\lambda w$ and the corresponding solution
$\mathbf{y}^\epsilon(g+ \lambda w)$
\begin{eqnarray*}
&&\int_D H^\epsilon(g+\lambda w) \sigma \left( \mathbf{y}^\epsilon(g+\lambda w)\right) 
: \nabla \mathbf{v}\,d\mathbf{x}
\nonumber\\
&&=\int_D H^\epsilon(g+\lambda w) \mathbf{f}\cdot \mathbf{v}\,d\mathbf{x}
+\int_{\Gamma_N} \mathbf{h}\cdot \mathbf{v}\,ds,\quad \forall \mathbf{v} \in W.
\end{eqnarray*}
Subtracting (\ref{stateeps}) from the above equation, we obtain
\begin{eqnarray*}
&&\int_D H^\epsilon(g+\lambda w) \sigma \left( \mathbf{y}^\epsilon(g+\lambda w)\right) 
: \nabla \mathbf{v}\,d\mathbf{x}
-\int_D H^\epsilon(g) \sigma \left( \mathbf{y}^\epsilon(g)\right) 
: \nabla \mathbf{v}\,d\mathbf{x}\nonumber\\
&&=\int_D H^\epsilon(g+\lambda w) \mathbf{f}\cdot \mathbf{v}\,d\mathbf{x}
-\int_D H^\epsilon(g) \mathbf{f}\cdot \mathbf{v}\,d\mathbf{x},\quad \forall \mathbf{v} \in W.
\end{eqnarray*}
Now, subtracting and adding the term 
$\int_D H^\epsilon(g+\lambda w) \sigma \left( \mathbf{y}^\epsilon(g)\right) 
: \nabla \mathbf{v}\,d\mathbf{x}$ in the first line,
dividing by $\lambda$, 
setting $\mathbf{z}^\epsilon_\lambda=\frac{\mathbf{y}^\epsilon(g+\lambda w)-\mathbf{y}^\epsilon(g)}{\lambda}$,
we get
\begin{eqnarray}
&&\int_D H^\epsilon(g+\lambda w) \sigma \left( \mathbf{z}^\epsilon_\lambda\right) 
: \nabla \mathbf{v}\,d\mathbf{x}
=
-\int_D \frac{H^\epsilon(g+\lambda w) - H^\epsilon(g)}{\lambda}\sigma \left( \mathbf{y}^\epsilon(g)\right) 
: \nabla \mathbf{v}\,d\mathbf{x}
\nonumber\\
&&+\int_D \frac{H^\epsilon(g+\lambda w) - H^\epsilon(g)}{\lambda}\mathbf{f}\cdot \mathbf{v}\,d\mathbf{x},
\quad \forall \mathbf{v} \in W.
\label{z_lambda}
\end{eqnarray}

Since $H^\epsilon\in \mathcal{C}^2(\mathbb{R})$, for each $y_0,h\in \mathbb{R}$,
there exists $\xi\in \mathbb{R}$ such that
$$
H^\epsilon(y_0+h)=H^\epsilon(y_0)+ h (H^\epsilon)^\prime(y_0)
+\frac{h^2}{2}(H^\epsilon)^{\prime\prime}(\xi),
$$
then, for each $\mathbf{x}\in \overline{D}$, there exists $\xi_x\in \mathbb{R}$ such that
$$
\frac{H^\epsilon(g(\mathbf{x})+\lambda w(\mathbf{x})) - H^\epsilon(g(\mathbf{x}))}{\lambda}
=w(\mathbf{x})(H^\epsilon)^\prime(g(\mathbf{x}))
+ \frac{\lambda}{2}w^2(\mathbf{x})(H^\epsilon)^{\prime\prime}(\xi_x).
$$
But, by construction $(H^\epsilon)^{\prime\prime}$ is bounded in $\mathbb{R}$ ($\epsilon$ is fixed) and
$w\in \mathcal{C}(\overline{D})$ is bounded in the compact $\overline{D}$. It follows that,
for all $\delta \in (0,1)$, there exists $\lambda_1(\delta)>0$ such that
\begin{equation}\label{uniform1}
\left\|\frac{H^\epsilon(g+\lambda w) - H^\epsilon(g)}{\lambda} -(H^\epsilon)^\prime(g)\,w
\right\|_{\mathcal{C}(\overline{D})}
\leq \delta,\quad \forall |\lambda| < \lambda_1(\delta),\ \lambda\neq 0.
\end{equation}
Then $\frac{H^\epsilon(g+\lambda w) - H^\epsilon(g)}{\lambda}$
converges to $(H^\epsilon)^\prime(g)\,w$ in $\mathcal{C}(\overline{D})$,
for $\lambda\rightarrow 0$.
We get that
\begin{equation}\label{uniform2}
\left\|\frac{H^\epsilon(g+\lambda w) - H^\epsilon(g)}{\lambda}\right\|_{\mathcal{C}(\overline{D})}
\leq M,\quad \forall |\lambda| < \lambda_1(\delta),\ \lambda\neq 0
\end{equation}
where $M=M(\epsilon)$ is independent of $\lambda$, but  depends on $\epsilon$.

As in Proposition \ref{prop:1}, we can obtain the estimate for $\mathbf{z}^\epsilon_\lambda\in W$
\begin{eqnarray}
\left\| \mathbf{z}^\epsilon_\lambda\right\|_{1,D} &\leq &
\frac{C\,M}{c(\epsilon)}\left( 
\left\| \sigma \left(\mathbf{y}^\epsilon(g)\right)\right\|_{0,D} 
+\left\| \mathbf{f}\right\|_{0,D} 
\right)
\nonumber\\
& \leq & \frac{C\,M}{c(\epsilon)}\left( 
C_1\left\| \mathbf{y}^\epsilon(g)\right\|_{1,D} 
+\left\| \mathbf{f}\right\|_{0,D} 
\right)
\label{norme_zlambda}
\end{eqnarray}
where $C_1>0$ is independent of $\lambda$, $\epsilon$
such that $\left\| \sigma \left(\mathbf{v}\right)\right\|_{0,D}
\leq C_1 \left\| \mathbf{v}\right\|_{1,D}$, for all $\mathbf{v}\in W$.

Let $\widetilde{\mathbf{z}}\in W$ such that, on a subsequence $\mathbf{z}^\epsilon_\lambda$
converges weakly to $\widetilde{\mathbf{z}}$ in $W$ and strongly in $\left(L^2(D)\right)^2$.
For passing to the limit on a subsequence in (\ref{z_lambda}),
{\color{black}we use Lemma }
6.1 from \cite{HMT2016}:
if $a,a_n\in L^\infty(D)$, $\left\| a_n\right\|_{0,\infty,D} \leq M$, $a_n\rightarrow a$ almost everywhere
in $D$, $b_n\rightarrow b$ weakly in $L^2(D)$ and $h\in L^2(D)$, then
$$
\lim_{n\rightarrow \infty}\int_D a_n b_n h \,dx=\int_D a\, b\, h \,dx.
$$
We can apply this Lemma for $a_n=H^\epsilon(g+\lambda_n w)$,
$b_n=\sigma \left( \mathbf{z}^\epsilon_{\lambda_n}\right)$ and $h=\nabla \mathbf{v}$.
By passing to the limit on a subsequence in (\ref{z_lambda}) we get that $\widetilde{\mathbf{z}}$
is solution of (\ref{limit}). But, as in Proposition \ref{prop:1}, we can show that
the problem (\ref{limit}) has a unique solution, then $\widetilde{\mathbf{z}}=\mathbf{z}$ and
$\mathbf{z}^\epsilon_\lambda$ converges to $\mathbf{z}$ for $\lambda\rightarrow 0$ without taking
subsequence, weakly in $W$ and strongly in $\left(L^2(D)\right)^2$.

Now, we will prove that $\mathbf{z}^\epsilon_\lambda$ converges to $\mathbf{z}$ strongly in $W$.
Subtracting (\ref{limit}) from (\ref{z_lambda}), we get
\begin{eqnarray}
&&\int_D H^\epsilon(g+\lambda w) \sigma \left( \mathbf{z}^\epsilon_\lambda\right) 
: \nabla \mathbf{v}\,d\mathbf{x}
-\int_D H^\epsilon(g) \sigma \left( \mathbf{z}\right) : \nabla \mathbf{v}\,d\mathbf{x}
\nonumber\\  
&=&-\int_D \left(
\frac{H^\epsilon(g+\lambda w) - H^\epsilon(g)}{\lambda}
-(H^\epsilon)^\prime(g) w \right)
\sigma \left( \mathbf{y}^\epsilon(g)\right) 
: \nabla \mathbf{v}\,d\mathbf{x}
\nonumber\\
&&+\int_D\left(
\frac{H^\epsilon(g+\lambda w) - H^\epsilon(g)}{\lambda}
-(H^\epsilon)^\prime(g) w\right)
\mathbf{f}\cdot \mathbf{v}\,d\mathbf{x}.
\label{substraction1}
\end{eqnarray}
Subtracting and adding the term 
$\int_D H^\epsilon(g) \sigma \left( \mathbf{z}^\epsilon_\lambda(g)\right) 
: \nabla \mathbf{v}\,d\mathbf{x}$ in the first line,
transferring some terms at the right-hand side, we get 
\begin{eqnarray}
&&
-\int_D H^\epsilon(g) \sigma \left( \mathbf{z}\right) : \nabla \mathbf{v}\,d\mathbf{x}
+\int_D H^\epsilon(g) \sigma \left( \mathbf{z}^\epsilon_\lambda\right) : \nabla \mathbf{v}\,d\mathbf{x}
\nonumber\\
&=&-\int_D \left(
\frac{H^\epsilon(g+\lambda w) - H^\epsilon(g)}{\lambda}
-(H^\epsilon)^\prime(g) w \right)
\sigma \left( \mathbf{y}^\epsilon(g)\right) 
: \nabla \mathbf{v}\,d\mathbf{x}
\nonumber\\
&&+\int_D\left(
\frac{H^\epsilon(g+\lambda w) - H^\epsilon(g)}{\lambda}
-(H^\epsilon)^\prime(g) w\right)
\mathbf{f}\cdot \mathbf{v}\,d\mathbf{x}
\nonumber\\
&&-\int_D \left(H^\epsilon(g+\lambda w) - H^\epsilon(g)\right)
\sigma \left( \mathbf{z}^\epsilon_\lambda\right) 
: \nabla \mathbf{v}\,d\mathbf{x}.
\label{substraction2}
\end{eqnarray}

At the left-hand side of (\ref{substraction2}), we have
$a(\mathbf{z}^\epsilon_\lambda -\mathbf{z}, \mathbf{v})$. 
Taking into account (\ref{uniform1}) and (\ref{uniform2}), the right-hand side of
(\ref{substraction2}) can be estimated by
\begin{eqnarray*}
&&\delta \left(  \left\| \sigma \left(\mathbf{z}^\epsilon_\lambda\right)\right\|_{0,D}
+\left\|\mathbf{f}\right\|_{0,D}
\right)\left\|\mathbf{v}\right\|_{0,D}
+\lambda  M \left\| \sigma \left(\mathbf{z}^\epsilon_\lambda\right)\right\|_{0,D}
\left\|\mathbf{v}\right\|_{0,D}
\\
&\leq & \delta \left( (1+M) C_1 \left\| \mathbf{z}^\epsilon_\lambda\right\|_{1,D}
+\left\|\mathbf{f}\right\|_{0,D}\right)\left\|\mathbf{v}\right\|_{0,D}
\end{eqnarray*}
for all $|\lambda| < min\left(\delta,\lambda_1(\delta)\right)$, $\lambda\neq 0$,
where $C_1>0$ such that $\left\| \sigma \left(\mathbf{v}\right)\right\|_{0,D}
\leq C_1 \left\| \mathbf{v}\right\|_{1,D}$, for all $\mathbf{v}\in W$.

Finally, from (\ref{norme_zlambda}), (\ref{norme_yeps}) and for
$\mathbf{v}=\mathbf{z}^\epsilon_\lambda -\mathbf{z}$,
we obtain
$$
\frac{c(\epsilon)}{C}\left\|\mathbf{z}^\epsilon_\lambda -\mathbf{z}\right\|_{1,D}^2
\leq
a(\mathbf{z}^\epsilon_\lambda -\mathbf{z}, \mathbf{z}^\epsilon_\lambda -\mathbf{z})
\leq \delta\, C_3(\epsilon) \left( \left\| \mathbf{f}\right\|_{0,D} +\left\| \mathbf{h}\right\|_{0,\Gamma_N}\right)
\left\|\mathbf{z}^\epsilon_\lambda -\mathbf{z}\right\|_{1,D}
$$
and after simplification, we get that
$\left\|\mathbf{z}^\epsilon_\lambda -\mathbf{z}\right\|_{1,D}
\leq \delta\frac{C\,C_3(\epsilon)}{c(\epsilon)}
\left( \left\| \mathbf{f}\right\|_{0,D} +\left\| \mathbf{h}\right\|_{0,\Gamma_N}\right)$
for all $\delta\in (0,1)$. 
Then $\mathbf{z}^\epsilon_\lambda$ converges to $\mathbf{z}$ strongly in $W$,
when $\lambda$ {\color{black}tends }
to $0$, but $\epsilon$ is fixed.

{\color{black}The linearity and the continuous dependence of its solution on the right-hand side
(on $w$) in equation (\ref{limit}) shows the G\^ateaux differentiability and ends the proof. }
\quad$\Box$

\begin{proposition}\label{prop:3}
 The {\color{black}directional derivative }
  of the objective function (\ref{jeps}) has the form
\begin{eqnarray}
 J^\prime(g)w
&=&
\int_D H^\epsilon(g) \mathbf{f}\cdot \mathbf{z}\,d\mathbf{x}
+\int_D (H^\epsilon)^\prime(g)\,w\, \mathbf{f}\cdot \mathbf{y}^\epsilon(g)\,d\mathbf{x}\nonumber\\
&&
+\int_{\Gamma_N} \mathbf{h}\cdot \mathbf{z}\,ds
+\ell \int_D  (H^\epsilon)^\prime(g)\,w\,d\mathbf{x}
\label{dj1}
\end{eqnarray}
for any $g$, $w$ in $X(D)$.
\end{proposition}

\noindent
\textbf{Proof.}
Let $g$, $w$ be fixed in $X(D)$ and $\lambda \neq 0$. We get
\begin{eqnarray}
\frac{J(g+\lambda w)-J(g)}{\lambda}
&=&
\frac{1}{\lambda}
\int_D \left( H^\epsilon(g+\lambda w) \mathbf{f}\cdot \mathbf{y}^\epsilon(g+\lambda w)
- H^\epsilon(g) \mathbf{f}\cdot \mathbf{y}^\epsilon(g)\right)d\mathbf{x}
\nonumber\\
&&+\int_{\Gamma_N} \mathbf{h}\cdot \frac{\mathbf{y}^\epsilon(g+\lambda w)-\mathbf{y}^\epsilon(g)}{\lambda} \,ds
\nonumber\\
&&+\ell \int_D \frac{ H^\epsilon(g+\lambda w)-H^\epsilon(g)}{\lambda} \,d\mathbf{x}
\label{gat_j}
\end{eqnarray}

From the Proposition \ref{prop:2}, 
$\mathbf{z}^\epsilon_\lambda=\frac{\mathbf{y}^\epsilon(g+\lambda w)-\mathbf{y}^\epsilon(g)}{\lambda}$
converges strongly to $\mathbf{z}$ in $W$. {\color{black}By the trace theorem on $\Gamma_N$,
we get }
$\mathbf{z}^\epsilon_\lambda |_{\Gamma_N}$ converges strongly to $\mathbf{z}|_{\Gamma_N}$ in 
$\left(L^2(\Gamma_N)\right)^2$. Consequently, the term of the second line in (\ref{gat_j})
converges to $\int_{\Gamma_N} \mathbf{h}\cdot \mathbf{z}\,ds$.

From (\ref{uniform1}), $\frac{H^\epsilon(g+\lambda w) - H^\epsilon(g)}{\lambda}$
converges uniformly to $(H^\epsilon)^\prime(g)\,w$ in $\mathcal{C}(\overline{D})$,
for $\lambda\rightarrow 0$, consequently, the term of the third line in (\ref{gat_j})
converges to $\ell \int_D  (H^\epsilon)^\prime(g)\,w\,d\mathbf{x}$.

It remains to study the right-hand side of the first line of (\ref{gat_j}).
Subtracting and adding $\int_D H^\epsilon(g+\lambda w) \mathbf{y}^\epsilon(g)
: \nabla \mathbf{v}\,d\mathbf{x}$, dividing by $\lambda$, we get
\begin{eqnarray*}
&&
\frac{1}{\lambda}\int_D \left( H^\epsilon(g+\lambda w)  \mathbf{y}^\epsilon(g+\lambda w)
- H^\epsilon(g+\lambda w) \mathbf{y}^\epsilon(g)\right)\cdot\mathbf{f}\, d\mathbf{x}\\
&&
+\frac{1}{\lambda}\int_D \left( H^\epsilon(g+\lambda w) \mathbf{y}^\epsilon(g)
- H^\epsilon(g) \mathbf{y}^\epsilon(g)\right) \cdot \mathbf{f}\, d\mathbf{x}\\
&=&
\int_D H^\epsilon(g+\lambda w) \frac{\mathbf{y}^\epsilon(g+\lambda w)- \mathbf{y}^\epsilon(g)}{\lambda} 
\cdot\mathbf{f}\, d\mathbf{x}\\
&& 
+\int_D \frac{H^\epsilon(g+\lambda w) - H^\epsilon(g)}{\lambda}
\mathbf{y}^\epsilon(g) \cdot \mathbf{f}\, d\mathbf{x}
\end{eqnarray*}
We have that $H^\epsilon(g+\lambda w)$ converges uniformly to $H^\epsilon(g)$ in $\mathcal{C}(\overline{D})$.
Using once again that $\mathbf{z}^\epsilon_\lambda$ converges strongly to $\mathbf{z}$ in $W$ and
$\frac{H^\epsilon(g+\lambda w) - H^\epsilon(g)}{\lambda}$ converges uniformly to $(H^\epsilon)^\prime(g)\,w$ 
in $\mathcal{C}(\overline{D})$, we get that the right-hand side of the first line of (\ref{gat_j})
converges to 
$\int_D H^\epsilon(g) \mathbf{z}\cdot \mathbf{f} \,d\mathbf{x}
+\int_D (H^\epsilon)^\prime(g)\,w\, \mathbf{y}^\epsilon(g)\cdot \mathbf{f}\,d\mathbf{x}$.
\quad$\Box$

We can give an expression of the {\color{black}directional derivative }
  of $J(g)$ without using $\mathbf{z}$.

\begin{proposition}\label{prop:4}
For any $g$, $w$ in $X(D)$, we have
\begin{eqnarray}
J^\prime(g)w
&=&
\int_D (H^\epsilon)^\prime(g) w \left[ 2\mathbf{f}\cdot \mathbf{y}^\epsilon(g) + \ell
-\sigma \left( \mathbf{y}^\epsilon(g)\right) : \nabla \mathbf{y}^\epsilon(g)
\right] d\mathbf{x}.
\label{dj2}
\end{eqnarray}
\end{proposition}

\noindent
\textbf{Proof.}
From (\ref{stateeps}), we put $\mathbf{v}=\mathbf{z}\in W$ and using
  {\color{black}(see (\ref{coercive})) }
$$
\sigma \left( \mathbf{y}^\epsilon(g)\right) : \nabla \mathbf{z}
=\lambda^S\left( \nabla \cdot \mathbf{y}^\epsilon(g)\right)\left(\nabla \cdot \mathbf{z}\right)
+2\mu^S \mathbf{e}\left( \mathbf{y}^\epsilon(g)\right)\mathbf{e}\left( \mathbf{z}\right)
=\sigma \left( \mathbf{z}\right) : \nabla \mathbf{y}^\epsilon(g)
$$
we get
\begin{eqnarray}
\int_D H^\epsilon(g)\sigma \left( \mathbf{z}\right) : \nabla \mathbf{y}^\epsilon(g)d\mathbf{x}
&=&
\int_D H^\epsilon(g)\, \mathbf{f}\cdot \mathbf{z}\, d\mathbf{x}
+\int_{\Gamma_N}\mathbf{h}\cdot \mathbf{z}\, ds.
\label{sigmaz}
\end{eqnarray}
Putting $\mathbf{v}=\mathbf{y}^\epsilon(g)$ in (\ref{limit}), it follows
\begin{eqnarray}
&&\int_D H^\epsilon(g) \sigma \left( \mathbf{z}\right) : \nabla \mathbf{y}^\epsilon(g)\,d\mathbf{x}
\nonumber\\
&=&
-\int_D (H^\epsilon)^\prime(g) w\, \sigma \left( \mathbf{y}^\epsilon(g)\right) : \nabla \mathbf{y}^\epsilon(g)\,d\mathbf{x}
+\int_D  (H^\epsilon)^\prime(g) w\, \mathbf{f}\cdot\mathbf{y}^\epsilon(g) \,d\mathbf{x}.
\label{limity}
\end{eqnarray}
From (\ref{sigmaz}) and (\ref{limity}), we obtain
\begin{eqnarray*}
&&\int_D H^\epsilon(g)\, \mathbf{f}\cdot \mathbf{z}\, d\mathbf{x}
+\int_{\Gamma_N}\mathbf{h}\cdot \mathbf{z}\, ds\\
&=&-\int_D (H^\epsilon)^\prime(g) w\, \sigma \left( \mathbf{y}^\epsilon(g)\right) : \nabla \mathbf{y}^\epsilon(g)\,d\mathbf{x}
+\int_D  (H^\epsilon)^\prime(g) w\, \mathbf{f}\cdot\mathbf{y}^\epsilon(g) \,d\mathbf{x}
\end{eqnarray*}
and taking into account (\ref{dj1}), we get the conclusion.
\quad$\Box$

\begin{remark} 
{\color{black}By the above result, one can obtain the form of the gradient of the cost and
avoid the use of an adjoint system. }
\end{remark}

In the following, we present some descent directions for the objective function.
We set
$$
d=2\mathbf{f}\cdot \mathbf{y}^\epsilon(g) + \ell
-\sigma \left( \mathbf{y}^\epsilon(g)\right) : \nabla \mathbf{y}^\epsilon(g)
$$
and since
$$
\sigma \left( \mathbf{y}^\epsilon(g)\right) : \nabla \mathbf{y}^\epsilon(g)
=\sigma \left( \mathbf{y}^\epsilon(g)\right) : \mathbf{e}(\mathbf{y}^\epsilon(g))
=\lambda^S(\nabla \cdot \mathbf{y}^\epsilon(g) )^2
+2\mu^S\mathbf{e}(\mathbf{y}^\epsilon(g)):\mathbf{e}(\mathbf{y}^\epsilon(g))
$$
then
\begin{equation}\label{d}
  d=2\mathbf{f}\cdot \mathbf{y}^\epsilon(g) + \ell
  -\left(
  \lambda^S(\nabla \cdot \mathbf{y}^\epsilon(g) )^2
+2\mu^S\mathbf{e}(\mathbf{y}^\epsilon(g)):\mathbf{e}(\mathbf{y}^\epsilon(g))
  \right).
\end{equation}
We have only $d\in L^1(D)$, in general.

We use, as in \cite{Tiba2018} the function $R:\mathbb{R}\rightarrow \mathbb{R}$
defined by
\begin{equation}\label{R}
R(r)=
\left\{
\begin{array}{ll}
c(-1+e^{r}), & r <0,\\
c(1-e^{-r}), & r \geq 0
\end{array}
\right.
\end{equation}
where $c>0$.
The function $R$ is strictly increasing, $R\left(\mathbb{R}\right) =]-c,c[$, 
$R(-r)=-R(r)$ and $r\,R(r)\geq 0$ for all $r\in \mathbb{R}$. 

\begin{proposition}\label{prop:5}
For $d$ given by (\ref{d}), {\color{black}the following are descent }
 directions for the objective function $J(g)$:
\begin{eqnarray}
& i)  & w_d=-H^\epsilon(g) d \label{3.i}\\
& ii) & w_d=-H^\epsilon(g) R(d) \label{3.ii}\\
& iii)& w_d=-\widetilde{d} \label{3.iii}
\end{eqnarray}
under the assumption that $w_d\in X(D)$.
At  $iii)$, $\widetilde{d}\in H^1(D)$ is the solution of
\begin{equation}\label{tilded}
\int_D  \gamma (\nabla \widetilde{d}\cdot \nabla v) + \widetilde{d}\,v\, d\mathbf{x}
=\int_D (H^\epsilon)^\prime(g) d\,v\, d\mathbf{x},
\quad \forall v\in H^1(D)
\end{equation}
and $\gamma>0$ is a parameter.
\end{proposition}

\noindent
\textbf{Proof.}
It is a consequence of Proposition \ref{prop:4}.
In the case i), we have
$$
J^\prime(g)w_d
=\int_D (H^\epsilon)^\prime(g) w_d d\, d\mathbf{x}
=-\int_D (H^\epsilon)^\prime(g) H^\epsilon(g) d^2 \, d\mathbf{x} < 0
$$
since $0 < (H^\epsilon)^\prime(r)$, $0 < H^\epsilon(r)$, for all $r\in \mathbb{R}$.

In the case ii), we have
$$
J^\prime(g)w_d
=-\int_D (H^\epsilon)^\prime(g) H^\epsilon(g)\, d\,R(d) \, d\mathbf{x} < 0
$$
since $r\,R(r) > 0$ for all $r\in \mathbb{R}^*$. 
In the case iii), we have
$$
J^\prime(g)w_d
=-\int_D (H^\epsilon)^\prime(g) d\,\widetilde{d} \, d\mathbf{x} 
=-\int_D  \gamma (\nabla \widetilde{d}\cdot \nabla \widetilde{d}) 
+ \widetilde{d}\,\widetilde{d}\, d\mathbf{x}
< 0
$$
since $\gamma>0$. 
\quad$\Box$

\begin{remark} 
More generally, $w_d=-\alpha\, d$ and $ w_d=-\alpha\, R(d)$,
where $\alpha \in L^\infty(D)$ and $\alpha \geq 0$, $\alpha \neq 0$, are descent directions, too.
{\color{black}For example, in the case $i)$, we have
  $J^\prime(g)w_d =-\int_D (H^\epsilon)^\prime(g) d^2 \alpha \, d\mathbf{x} <0$, since
$0 < (H^\epsilon)^\prime(r)$ for all $r\in \mathbb{R}^*$, $d^2\geq 0$, $\alpha \geq 0$. }
The case $iii)$ is inspired by \cite{Burger2003}.
\end{remark}

\section{Numerical examples}
\setcounter{equation}{0}

We have employed the software FreeFem++, \cite{freefem++}.
The dimensions and the starting domains are from the web site of the team
directed by G. Allaire \cite{Allaire}, the files \texttt{levelset-cantilever.edp} and 
\texttt{pont.homog.struct.edp}.\\
Our approach decreases the cost and ensures both boundary and/or topology variations, 
including the creation of new holes.

\bigskip

\textbf{Algorithm}

\vspace{2mm}

\textbf{Step 1} Let $g_0$ be the initial guess.
{\color{black}Fix a bound $N$ }
 for the number of iterations and put $n:=0$.

 \textbf{Step 2} Let $w_n$ be a descent direction of $J$ in $g_n$,
 given by {\color{black}(\ref{3.i}), (\ref{3.ii}), or (\ref{3.iii}). }

\textbf{Step 3} 
Compute the {\color{black}directional derivative }
$J^\prime(g_n)w_n$  according to Proposition \ref{prop:4}.\\
If $J^\prime(g_n)w_n=0$, then \textbf{Stop}.

\textbf{Step 4} Find 
$$
\lambda_n\in \arg \min_{\lambda\in \mathbb{R}} J(g_n+\lambda\,w_n)
$$
obtained via some line search and put $g_{n+1}=g_n+\lambda_n\,w_n$.
{\color{black} Practically, we look for $\lambda = \rho^i$,
$\rho\in (0,1)$ and $i=0,1,\dots$. The maximal number of iterations for the line search is fixed
  to 10. Alternatively, we can use backtracking line search method,
  see \cite{Dennis1996} or other method.
}

\textbf{Step 5} If $n+1=N$ then \textbf{Stop}

\textbf{Step 6} If $| J(g_{n})-J(g_{n+1}) | < tol$ then \textbf{Stop}, \\
else update $n:=n+1$ and go to \textbf{Step 2}.

\bigskip

{\color{black}The final domains given by the above algorithm  are  not necessary globally optimal.
In the following, optimal domain means final computed domain  in the descent procedure.}

\bigskip

In the previous sections, we have supposed that $\Gamma_D$ and  $\Gamma_N$ are given.
For the numerical tests, we assume  that only $\Gamma_N \subset \partial D$ is given, but
$\Gamma_D$ is unknown. The space $W$ for the weak formulation 
(\ref{stateeps}) is 
$$
W=\{\mathbf{v} \in \left(H^1(D)\right)^2;\ \mathbf{v}=0\hbox{ on }\Sigma_D \}.
$$
Generally, $\Sigma_D$ is not a subset of $\partial\Omega_g$ and the Dirichlet boundary condition
(\ref{elast2}) is imposed in fact only on $\partial\Omega_g\cap \Sigma_D$.
We set $\Gamma_D=\partial\Omega_g\cap \Sigma_D$. {\color{black}In terms of the control $g$, this can be
  ensured by imposing in the definition of $X(D)$ that
  $g(\mathbf{x}) < 0, \mathbf{x} \in \partial D \setminus [\Sigma_D \cup \Gamma_N] $ too
(compare Remark \ref{rem:3.1} ). }

\noindent
If $\Omega_g\subset\subset D$, then $\mathbf{x}\in \partial\Omega_g$ yields
$g(\mathbf{x})=0$. Otherwise, when 
$meas\left(\partial\Omega_{g}\cap \partial D\right)>0$,
it is possible $\mathbf{x}\in \left(\partial\Omega_g \cap \partial D\right)$ and $g(\mathbf{x}) >0$.

The initial parametrization here satisfies $g_0(\mathbf{x}) > 0,\ \mathbf{x}\in \Gamma_N$ and
$\Gamma_N\subset\partial\Omega_{g_0}$.
Also, we have $meas\left(\partial\Omega_{g_0}\cap \Sigma_D\right)>0$.
For particular initial parametrizations used here, we have observed that
$g_n(\mathbf{x}) > 0,\ \mathbf{x}\in \Gamma_N$ and $meas\left(\partial\Omega_{g_n}\cap \Sigma_D\right)>0$
for all the iterations $n$ until convergence.
For different initial parametrization, this property may not hold systematically and should be
imposed as a constraint on $g$.

\bigskip

\textbf{Example 1. Cantilever} 

\begin{figure}[ht]
  \begin{center}
\includegraphics[width=5.5cm]{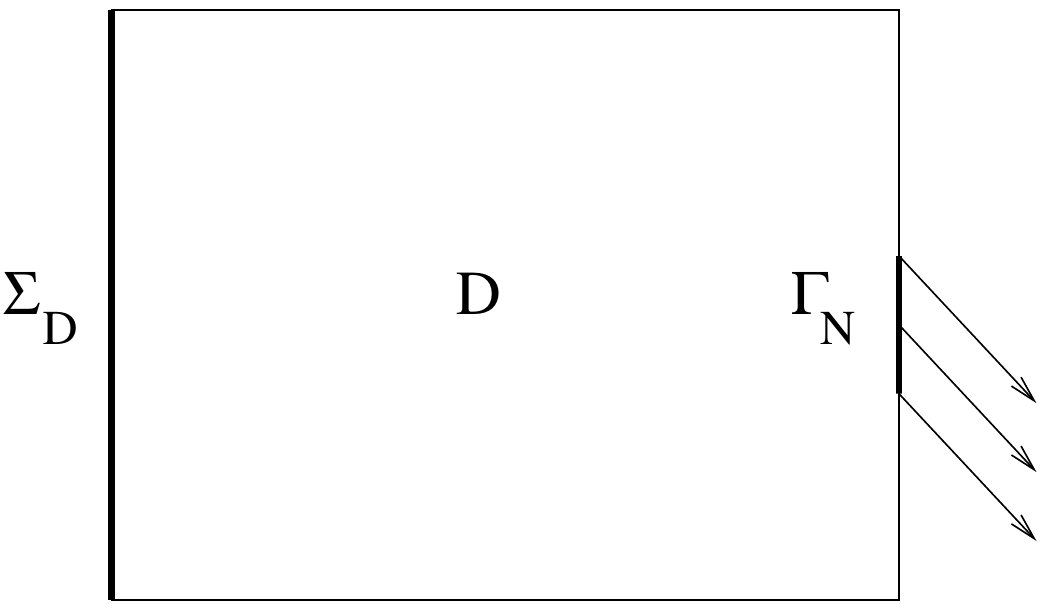}\quad    
\includegraphics[width=6cm]{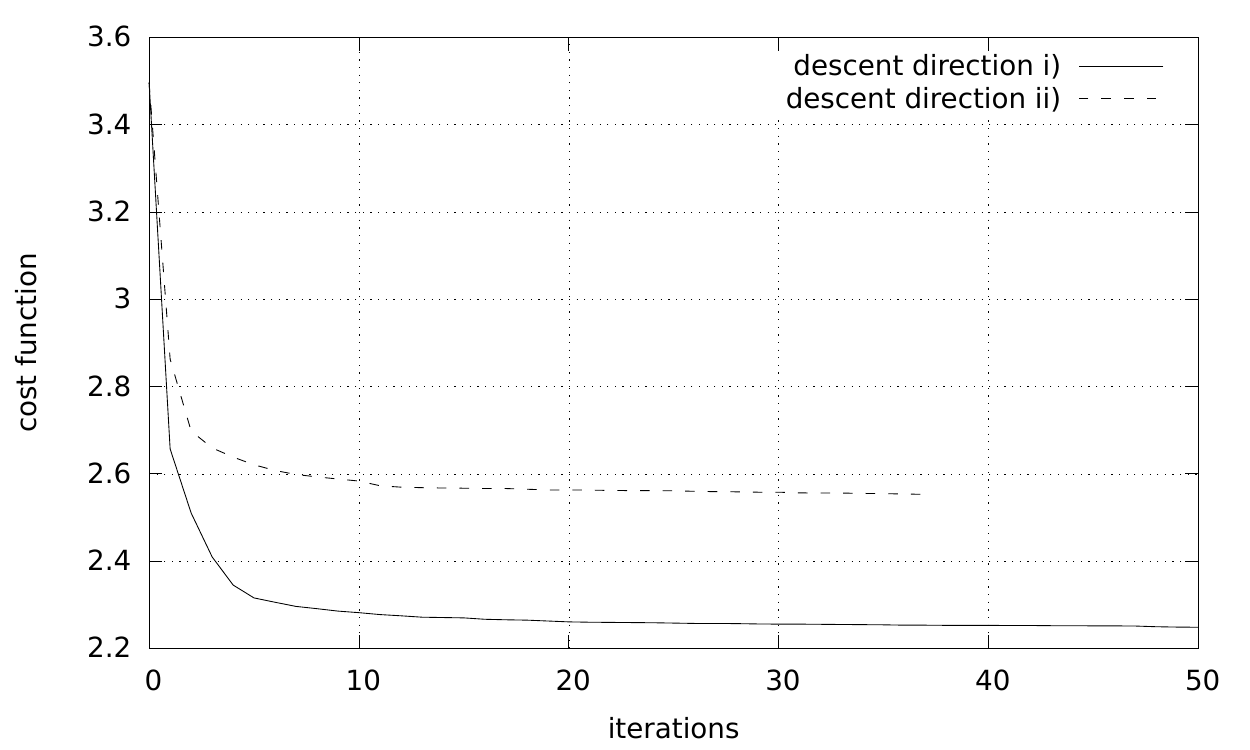}
\end{center}
  \caption{Cantilever. Left: Geometrical configuration of $D$.
    Right: Convergence history of the objective functions for descent directions
$i)$ given by (\ref{3.i}) and $ii)$ given by (\ref{3.ii}).}
\label{fig:canti_Jk}
\end{figure}

We have $D=]0,2[ \times ]-0.5,0.5[$, $\Sigma_D=\{0\} \times ]-0.5,0.5[$,
$\Gamma_N=\{2\} \times ]-0.1,0.1[$, see Figure \ref{fig:canti_Jk}, left.
We work with Lam\'e coefficients $\lambda^S=1$, $\mu^S=8$
{\color{black} and $\rho=0.6$. }
The volume load is $\mathbf{f}=(0,0)$, the surface load on $\Gamma_N$ is $\mathbf{h}=(0,-5)$
and the parameter in the objective function associated to the volume of the structure
is $\ell=0.5$.

We use for $D$ a mesh of 45638 triangles and 23120 vertices. 
For the approximation of $g$ and $H^\epsilon(g)$ we use piecewise linear 
finite element, globally continuous and for $\mathbf{y}$ we use 
the finite element $\mathbb{P}_2$, piecewise polynomial of degree two, see \cite{Raviart2004}.
We set $\epsilon=10^{-2}$ the penalization parameter and $tol=10^{-6}$, $N=50$ for the
stopping tests.

\begin{figure}[ht]
\begin{center}
\includegraphics[width=6cm]{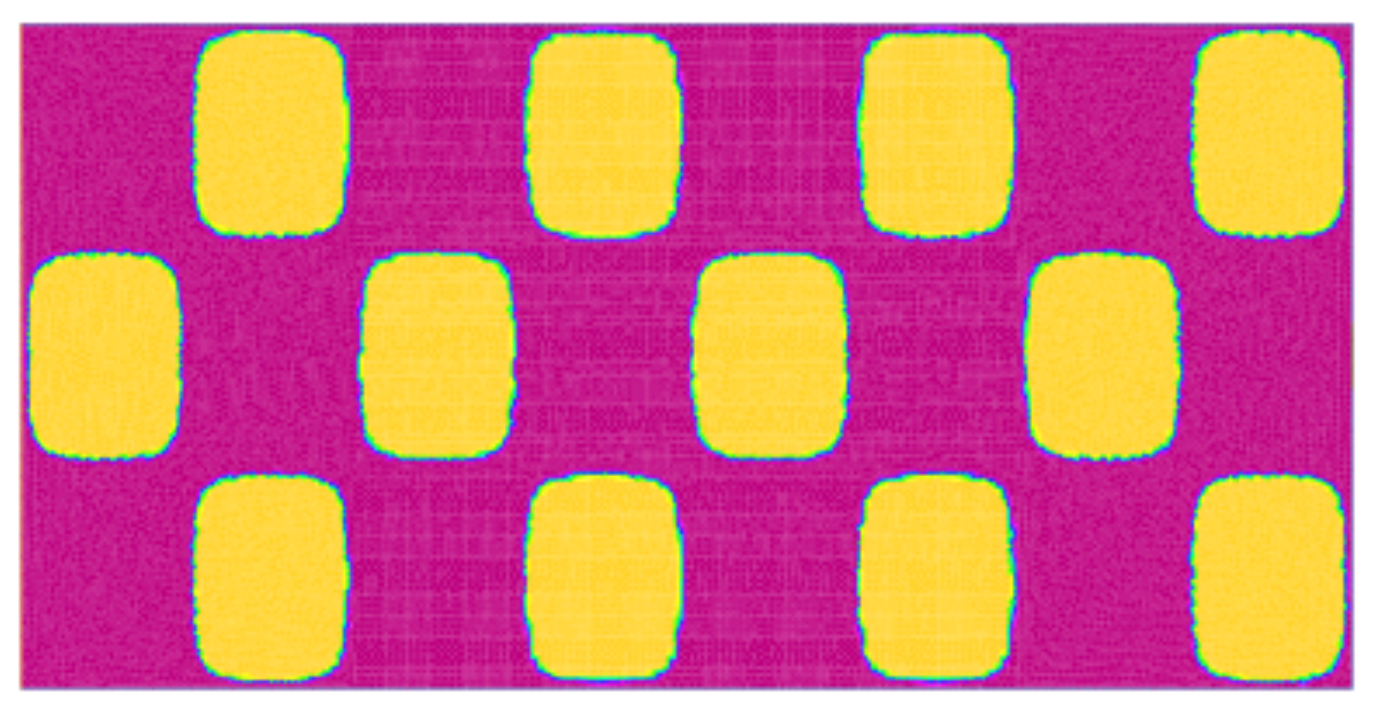}
\quad
\includegraphics[width=6cm]{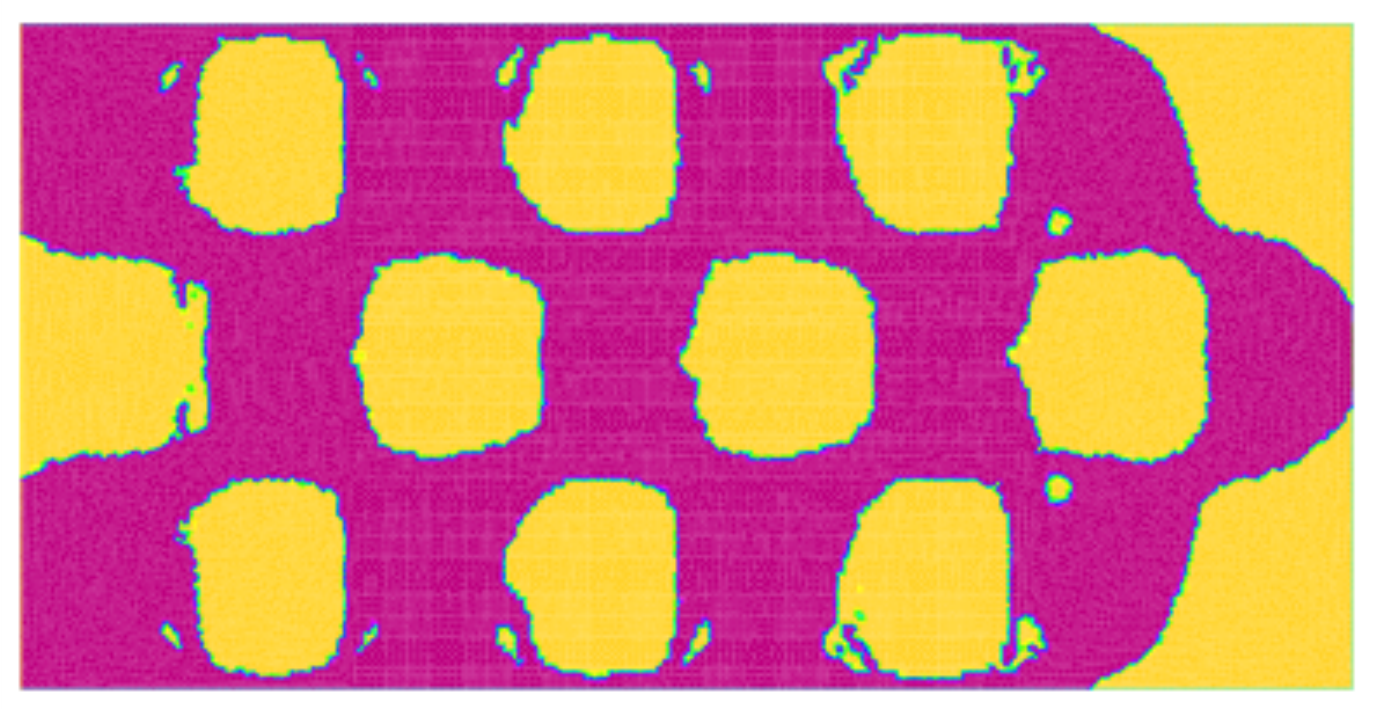}\\
\includegraphics[width=6cm]{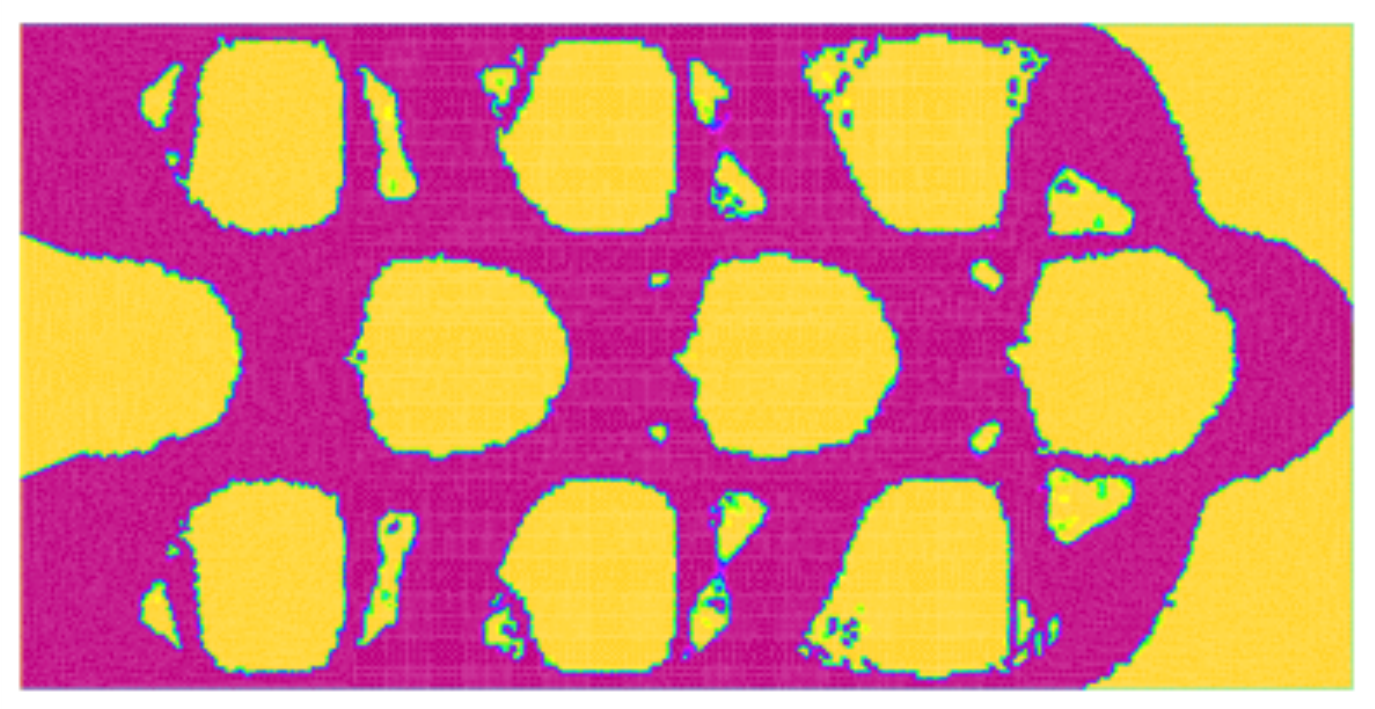}
\quad
\includegraphics[width=6cm]{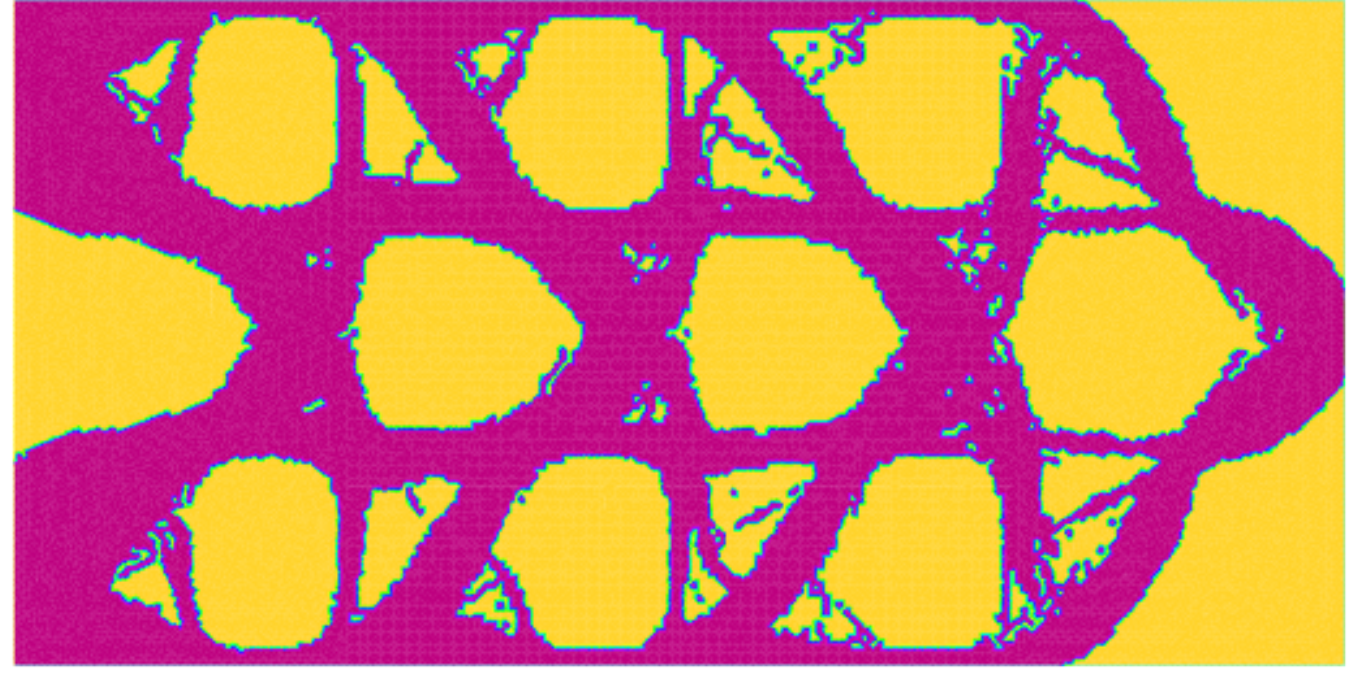}
\end{center}
\caption{Cantilever. Initial (top, left), intermediate and optimal (bottom, right, after 50 iterations)
  domains using descent direction $i)$.}
\label{fig:canti_init_final}
\end{figure}

\begin{figure}[ht]
\begin{center}
\includegraphics[width=6cm]{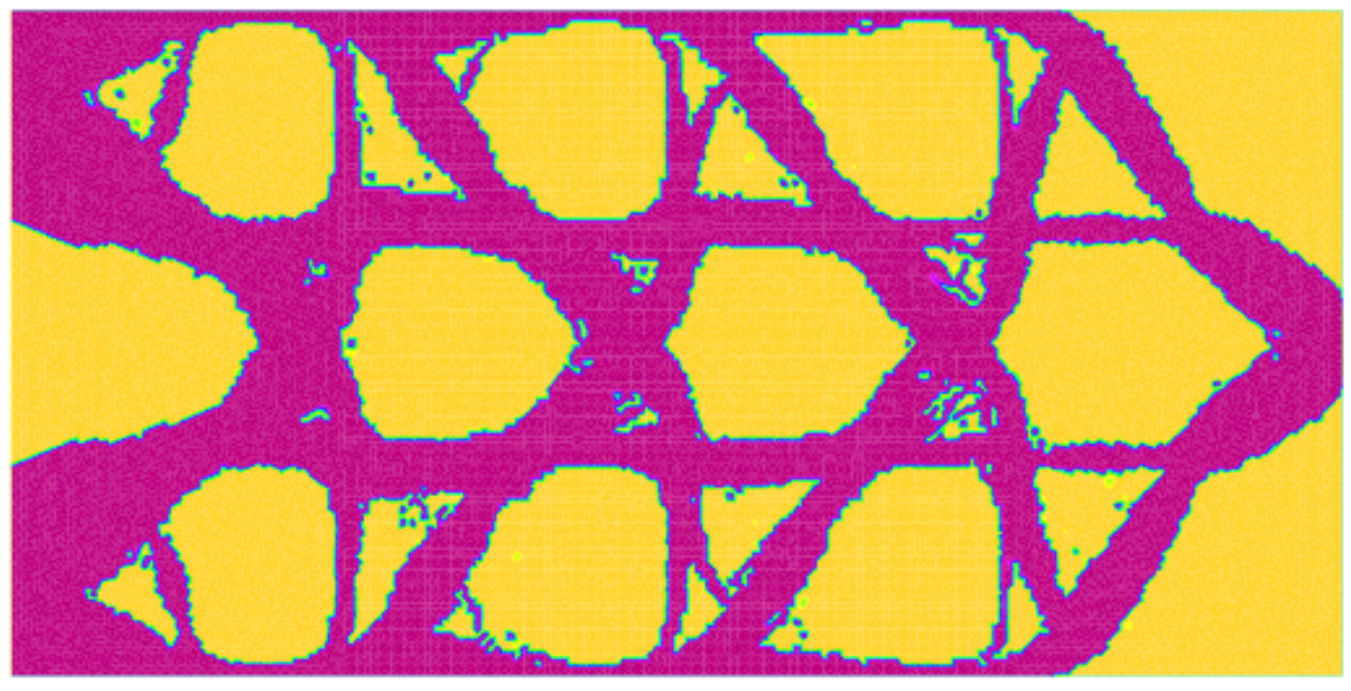}
\quad
\includegraphics[width=6cm]{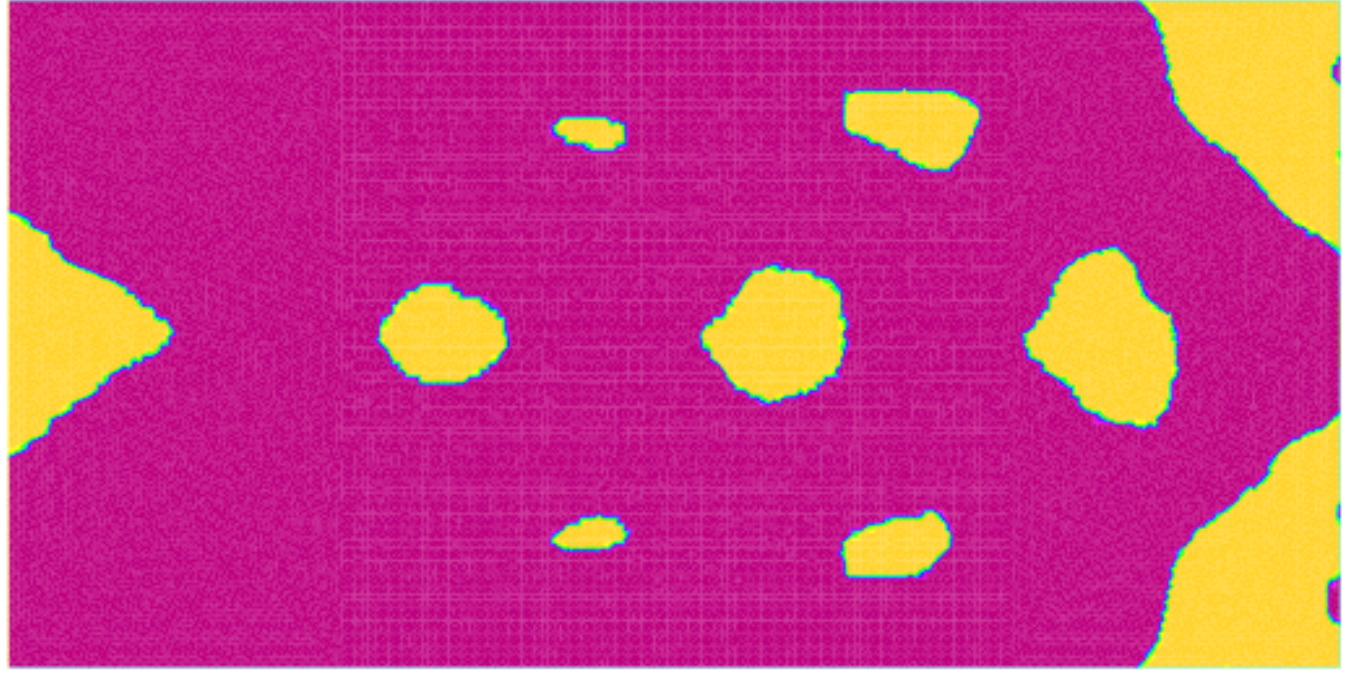}
\end{center}
\caption{Cantilever. Optimal domains using descent directions
  $ii)$ given by (\ref{3.ii}) (left, after 38 iterations) and $iii)$ given by (\ref{3.iii})
  (right, after 3 iterations)
  for the initial domain as in Figure \ref{fig:canti_init_final}.}
\label{fig:canti_R_reg}
\end{figure}

The initial domain is obtained for 
$g_0(x_1,x_2)=0.1- \sin(4\pi x_1)\sin\left( 3\pi(x_2-0.5)\right)$ and the initial
value of the objective function is $J(g_0)=3.49524$.
{\color{black} This starting domain is also used by  \cite{Allaire}. It has many initial holes
  and the algorithm ``closes'' some of them, but also produces new holes as may be seen in
  Figure \ref{fig:canti_init_final} and \ref{fig:canti_R_reg}. }

The history of the objective functions for descent directions
$i)$ and $ii)$ given by (\ref{3.i}) and (\ref{3.ii})
is presented in Figure \ref{fig:canti_Jk} right,
the optimal value is $2.24849$ in the case $i)$ after 50 iterations (\textbf{Step 5})
and $2.55336$ in the case $ii)$ after 38 iterations (\textbf{Step 6}).
The stopping test \textbf{Step 3} is obtained for $n=3$ for the descent direction $iii)$ 
given by (\ref{3.iii}), (\ref{tilded}) with $\gamma=0.001$, 
the values of the objective function are:
$J(g_0)=3.49524$, $J(g_1)=1.45725$, $J(g_2)=1.45704$, $J(g_3)=1.45626$.
The initial, intermediate and  optimal domains using different descent directions
are presented in Figure \ref{fig:canti_init_final} and \ref{fig:canti_R_reg}.
{\color{black}The final value of the objective function in the case $iii)$ ($1.45626$)
  is less that in the case $i)$ ($2.24849$) which is less than in the
  case $ii)$ ($2.55336$). We also observe that the volume of final domain in the case $iii)$ is
  larger than in the case $i)$ or  $ii)$.
  These results are influenced by $\ell=0.5$ which means that the material is ``cheap''.
  In the case where $\ell$ is large, in other words, the material is ``expensive'',
  the final domains will have smaller volumes. The acting forces are the same and smaller
  volumes allow larger displacements and larger compliances.
}

We have tested the dependence on $\epsilon$ of the optimal solution given by the algorithm presented
in this paper.
We use the same mesh, the same initial domain given by 
$g_0(x_1,x_2)=0.1- \sin(4\pi x_1)\sin\left( 3\pi(x_2-0.5)\right)$ and the descent direction
$i)$ given by (\ref{3.i}).

For small $\epsilon$, the value of $c(\epsilon)$ from Proposition \ref{prop:1} is close to zero, then
solving numerically the linear system associated to (\ref{stateeps}) is difficult.
We have replaced in (\ref{stateeps}), $H^\epsilon$ by
$$
\widehat{H}^\epsilon(r)=
\left\{
\begin{array}{ll}
1-\frac{1}{2}e^{-\frac{r}{\epsilon}}, & r \geq 0,\\
\max\left(0.0001,\frac{1}{2}e^{\frac{r}{\epsilon}}\right), & r < 0 .
\end{array}
\right.
$$

For $\epsilon=10^{-3}$, the initial objective function is $3.52187$.
The stopping test \textbf{Step 3} is obtained for
$n=38$, the final objective function is $2.29428$ and we denote the final displacement by
$\mathbf{y}^{0.001}$. For $\epsilon=10^{-4}$, the initial objective function is $3.54231$,
the stopping test \textbf{Step 3} is obtained for
$n=28$, the final objective function is $2.37167$. Similarly, we denote by $\mathbf{y}^{0.0001}$
the final displacement.
The final domains are similar to the case $\epsilon=10^{-2}$, Figure
\ref{fig:canti_init_final}, (bottom, right). The final displacement is denoted by
$\mathbf{y}^{0.01}$ when $\epsilon=10^{-2}$.
We have computed the differences of the final displacements in norms $L^2$ and $H^1$:\\
$\left\| \mathbf{y}^{0.01}-\mathbf{y}^{0.0001}\right\|_{L^2(\Omega_{28})}
=0.083864$,
$\left\| \mathbf{y}^{0.01}-\mathbf{y}^{0.0001}\right\|_{H^1(\Omega_{28})}
=0.382696$,\\
$\left\| \mathbf{y}^{0.001}-\mathbf{y}^{0.0001}\right\|_{L^2(\Omega_{28})}
=0.066597$,
$\left\| \mathbf{y}^{0.001}-\mathbf{y}^{0.0001}\right\|_{H^1(\Omega_{28})}
=0.362007$.


\bigskip

\textbf{Example 2. Bridge} 

\begin{figure}[ht]
  \begin{center}
\includegraphics[width=4cm]{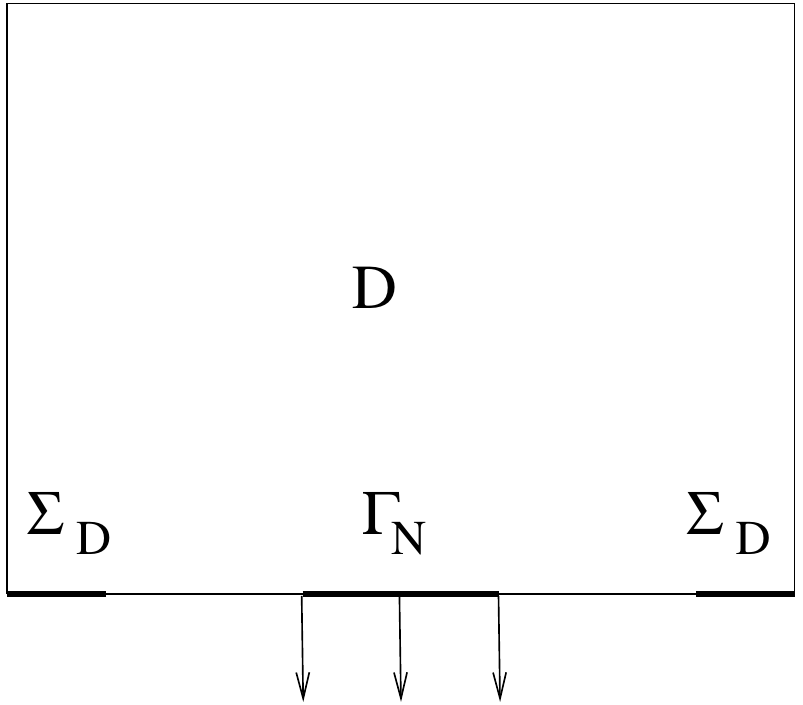}\quad     
\includegraphics[width=6cm]{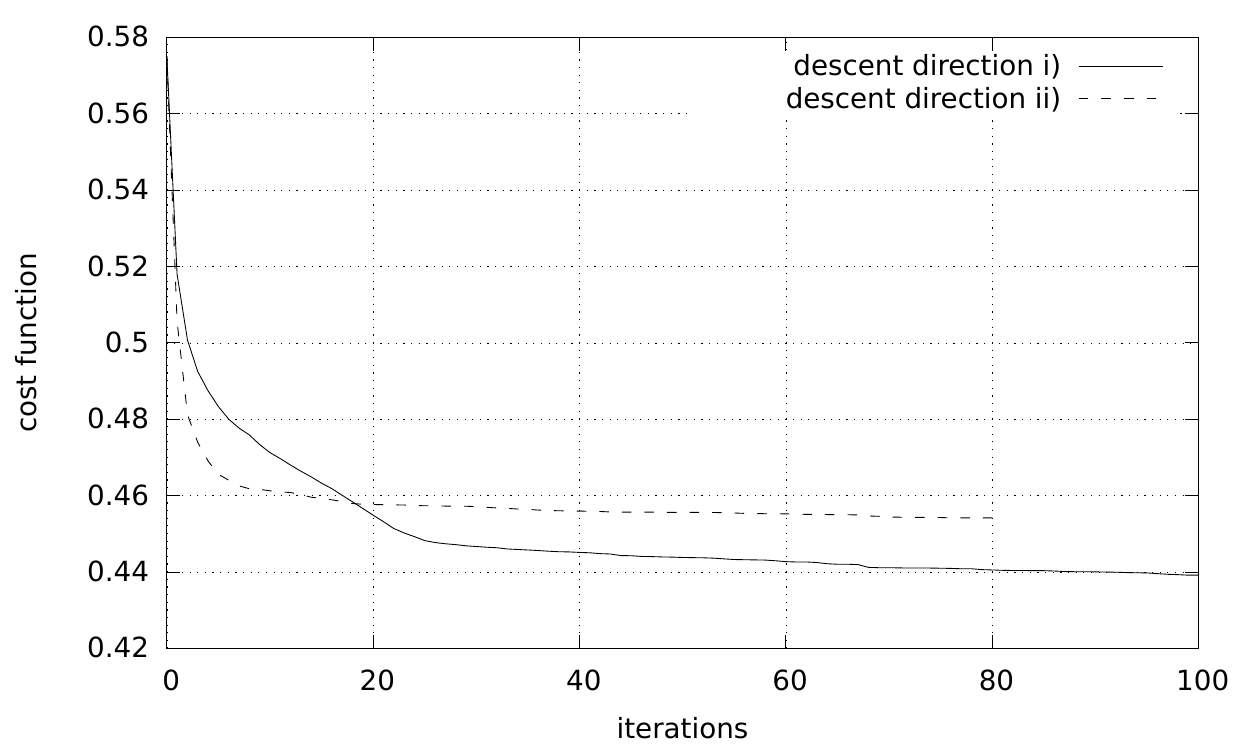}
\end{center}
  \caption{Bridge. Left: Geometrical configuration of $D$.
    Right: Convergence history of the objective functions for descent directions
$i)$ and $ii)$.}
\label{fig:pont2_Jk}
\end{figure}

We have $D=]-1,1[ \times ]0,1.2[$,
$\Sigma_D=\left( ]-1,-0.9[ \cup ]0.9,1[ \right) \times \{0\} $,\\
$\Gamma_N= ]-0.1,0.1[ \times \{0\} $, see Figure \ref{fig:pont2_Jk} left.
We work with Young modulus $E=1$, Poisson ratio $\nu=0.3$
{\color{black} and $\rho=0.6$. }
The volume load is $\mathbf{f}=(0,0)$, the surface load on $\Gamma_N$ is $\mathbf{h}=(0,-1)$
and the parameter in the objective function associated to the volume of the structure
is $\ell=0.1$.

We use for $D$ a mesh of 54510 triangles and 27576 vertices.
We set $\epsilon=10^{-2}$ {\color{black}for }
  the penalization parameter and $tol=10^{-6}$, $N=100$ for the
stopping tests.

\begin{figure}[ht]
\begin{center}
\includegraphics[width=6cm]{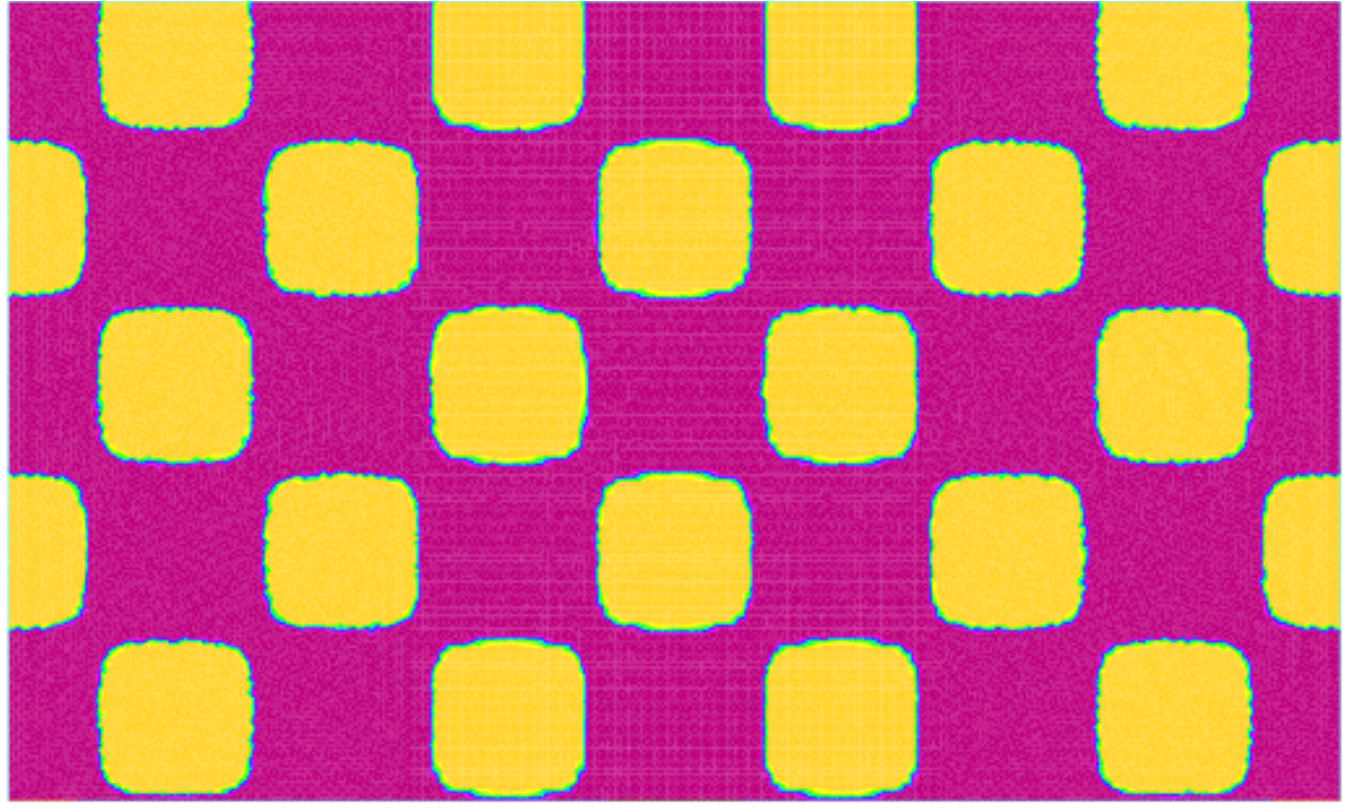}
\quad
\includegraphics[width=6cm]{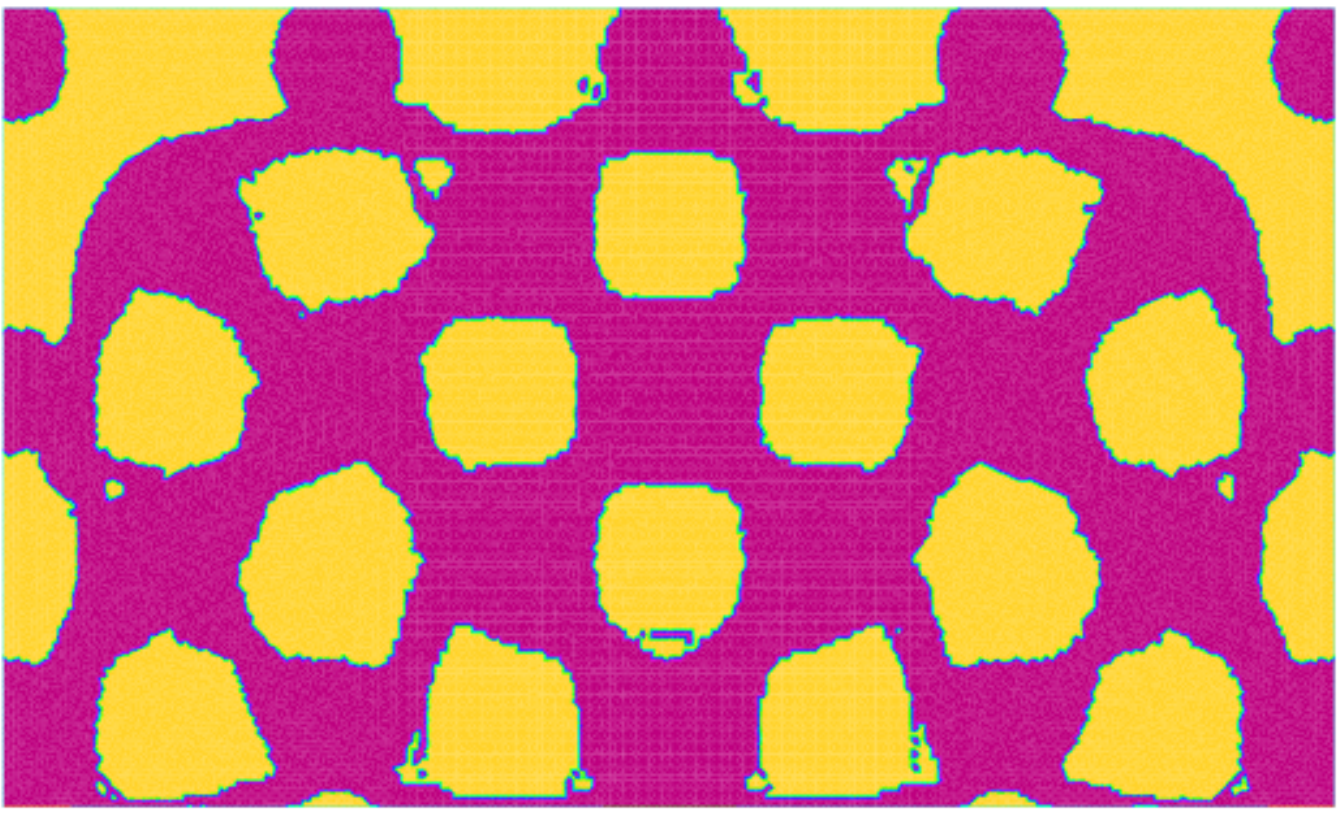}\\
\includegraphics[width=6cm]{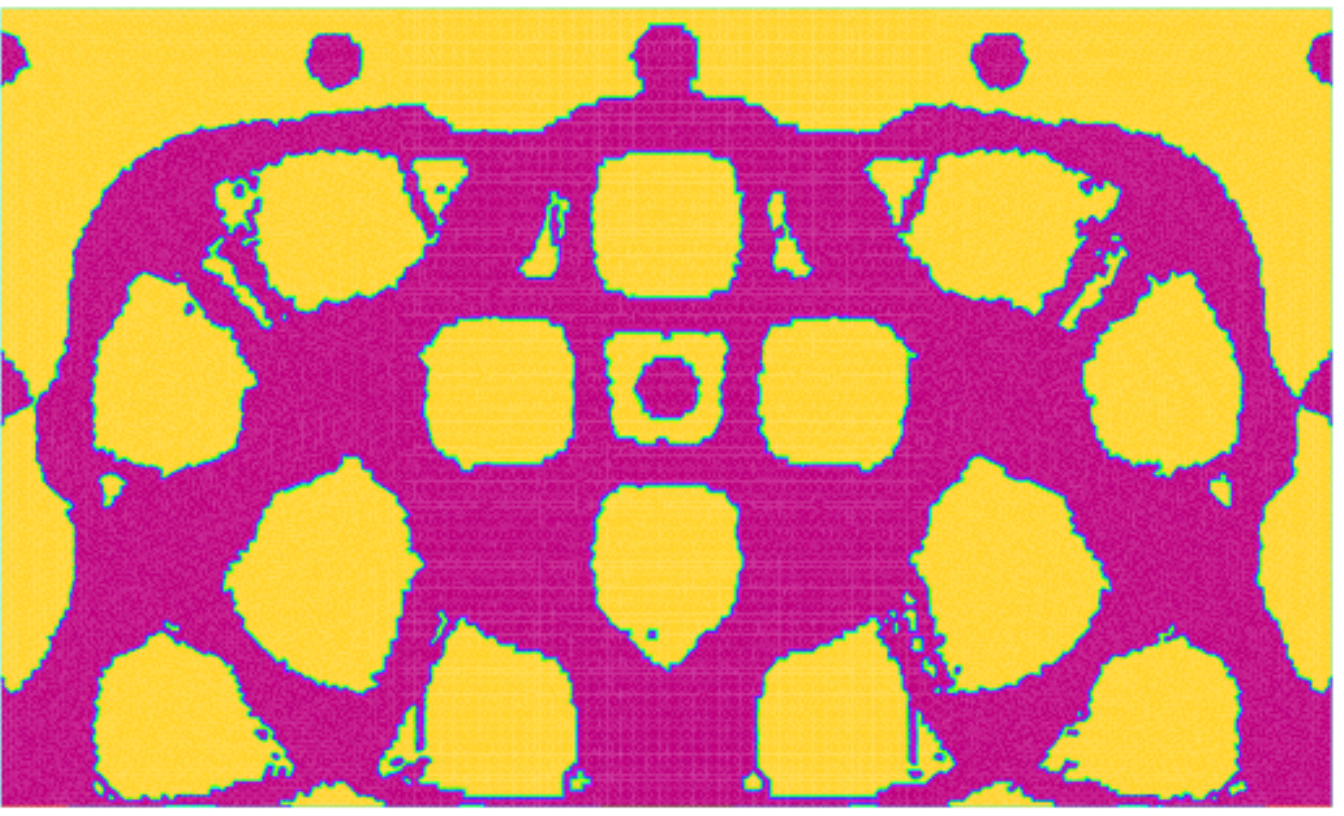}
\quad
\includegraphics[width=6cm]{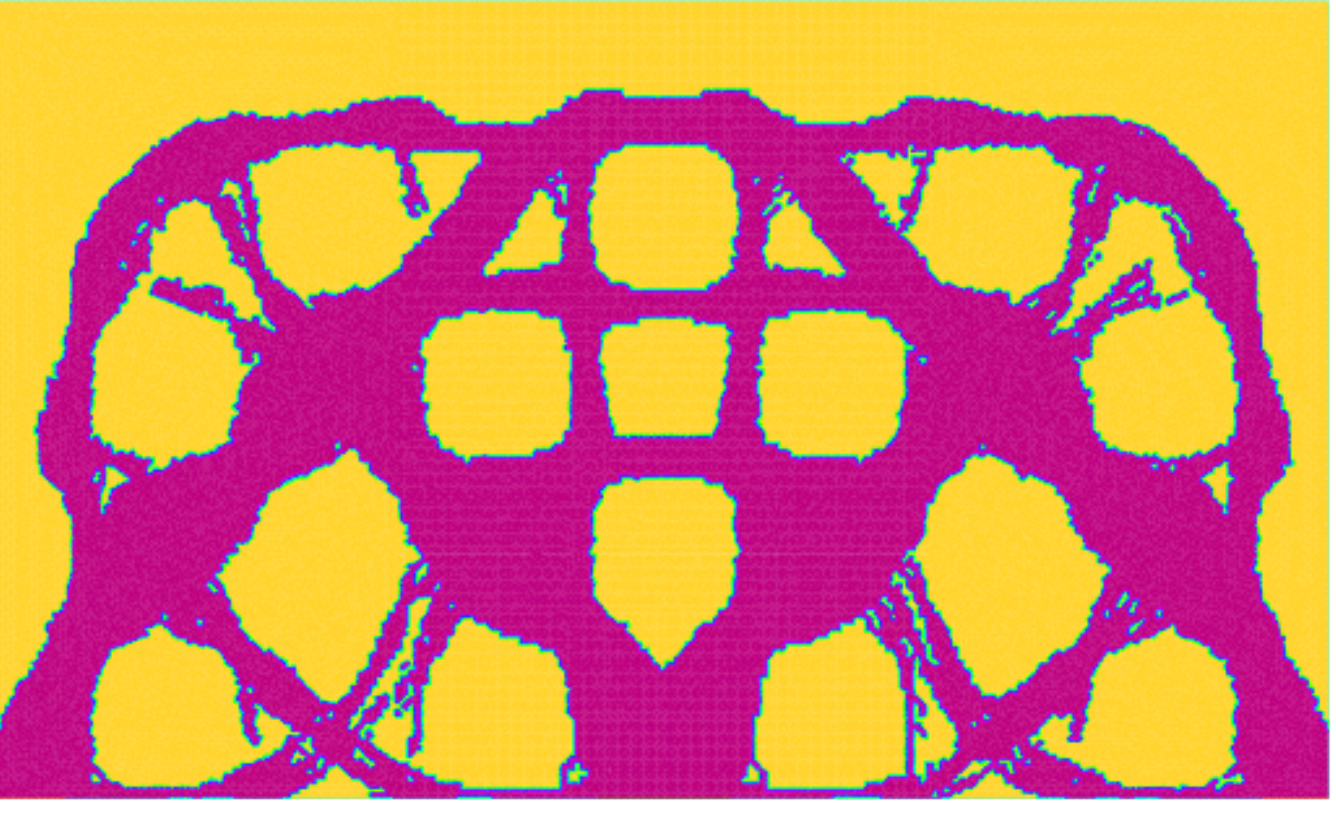}
\end{center}
\caption{Bridge. Initial (top, left), intermediate and optimal (bottom, right, after 100 iterations)
  domains using descent direction $i)$.}
\label{fig:pont2_init_final}
\end{figure}

\begin{figure}[ht]
\begin{center}
\includegraphics[width=6cm]{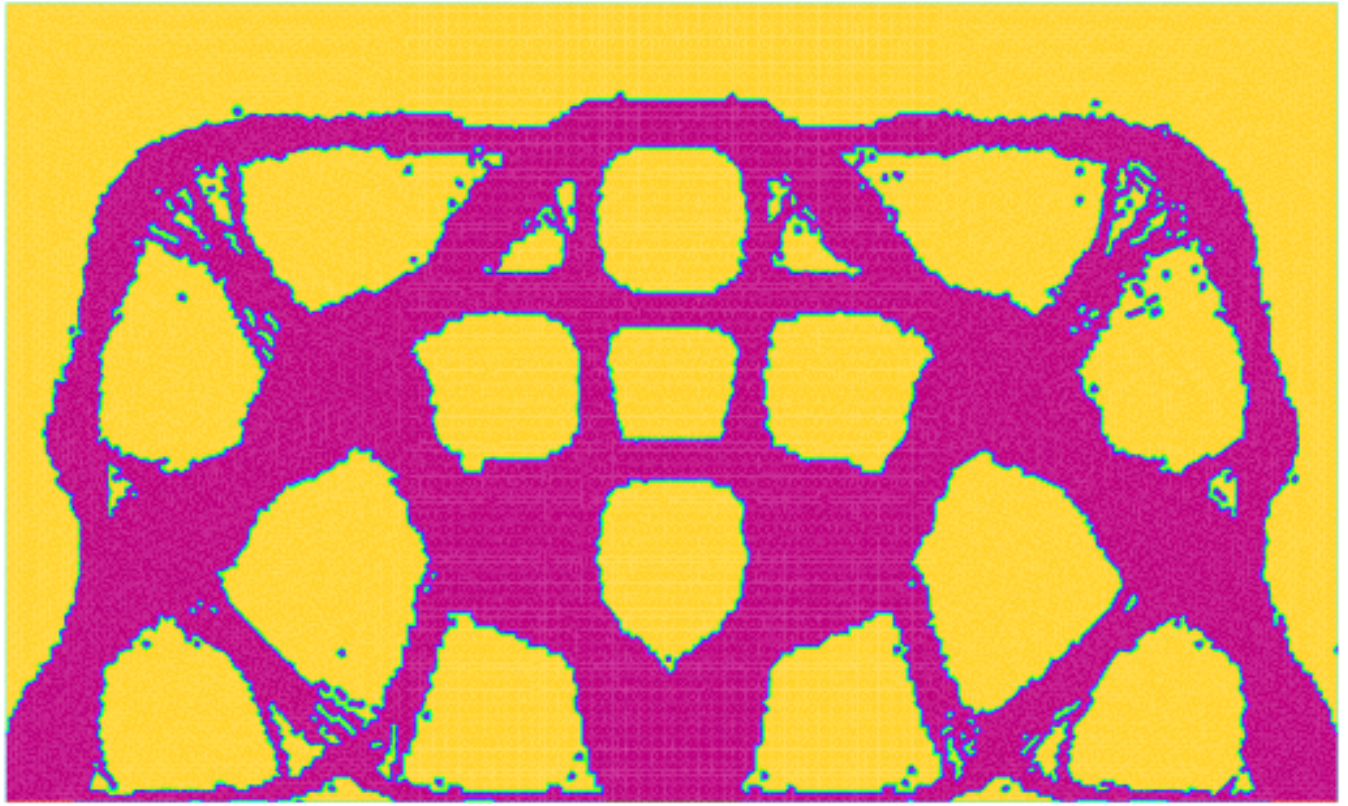}
\end{center}
\caption{Bridge. Optimal domain using descent directions
$ii)$ after 80 iterations, for the initial domain as in Figure \ref{fig:pont2_init_final}.}
\label{fig:pont2_sigma_R}
\end{figure}

The initial domain is obtained for 
$$
g_0(x_1,x_2)=0.1- \sin\left(4\pi (x_1-0.125)\right)\sin\left( 4\pi(x_2-0.5)\right)
$$
and the initial value of the objective function is $J(g_0)=0.574918$.
The algorithm stops after 100 iterations  (\textbf{Step 5}) when using descent directions
$i)$ and after 80 iterations  (\textbf{Step 6}) when using descent directions $ii)$,
the optimal value of the objective function
 is $0.43918$ in the case $i)$ and $0.454161$ in the case $ii)$.
The initial, intermediate and the optimal domains using different descent directions
are presented in Figure \ref{fig:pont2_init_final} and \ref{fig:pont2_sigma_R}.
The descent direction
$iii)$ did not work properly in this example.

\begin{figure}[ht]
\begin{center}
\includegraphics[width=6cm]{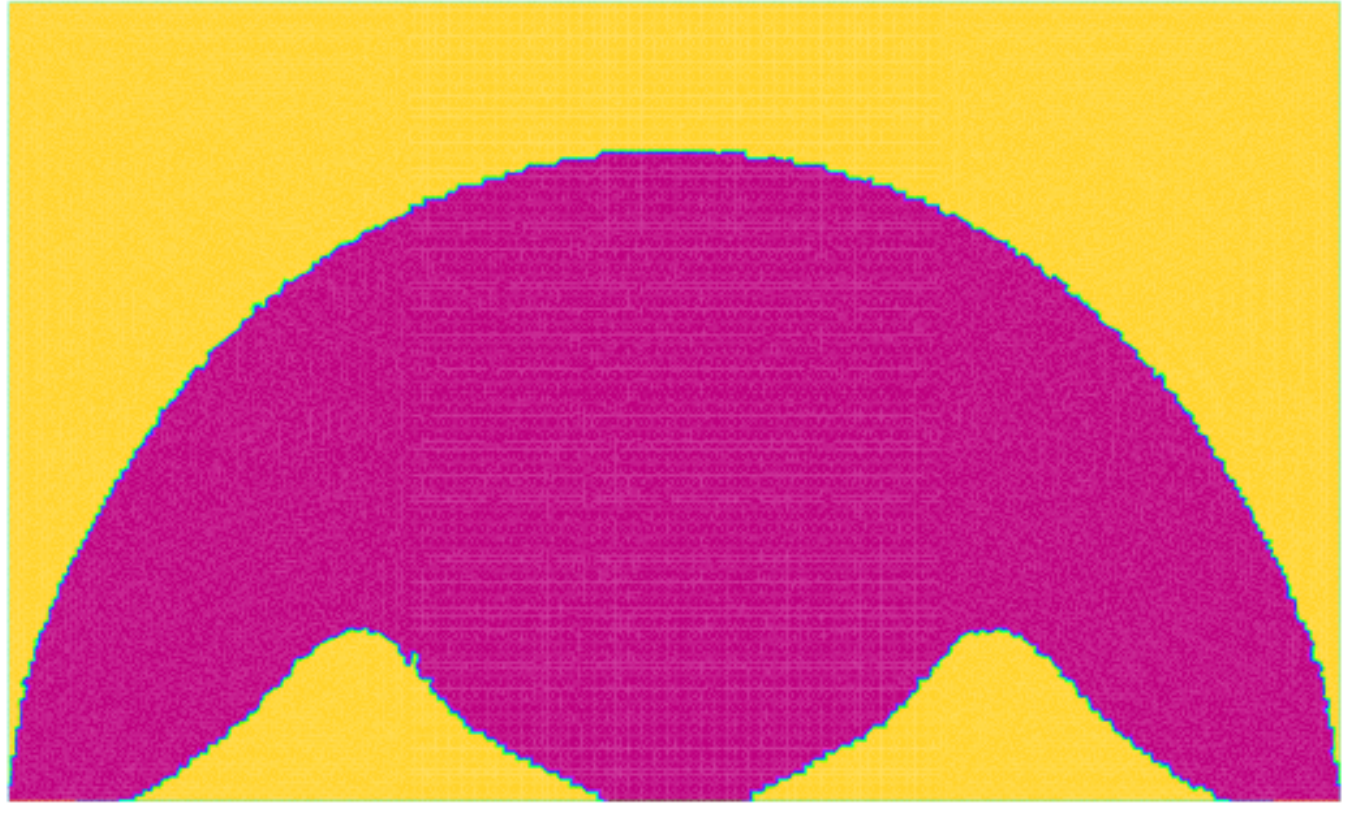}
\end{center}
\caption{Bridge. Optimal domain using descent directions
  $i)$ after 100 iterations, for the initial domain $\Omega_0=]-1,1[ \times ]0,0.6[$,
  the bottom half of $D$.}
\label{fig:pont2_sigma_half_D}
\end{figure}

\begin{figure}[ht]
  \begin{center}
\includegraphics[width=6cm]{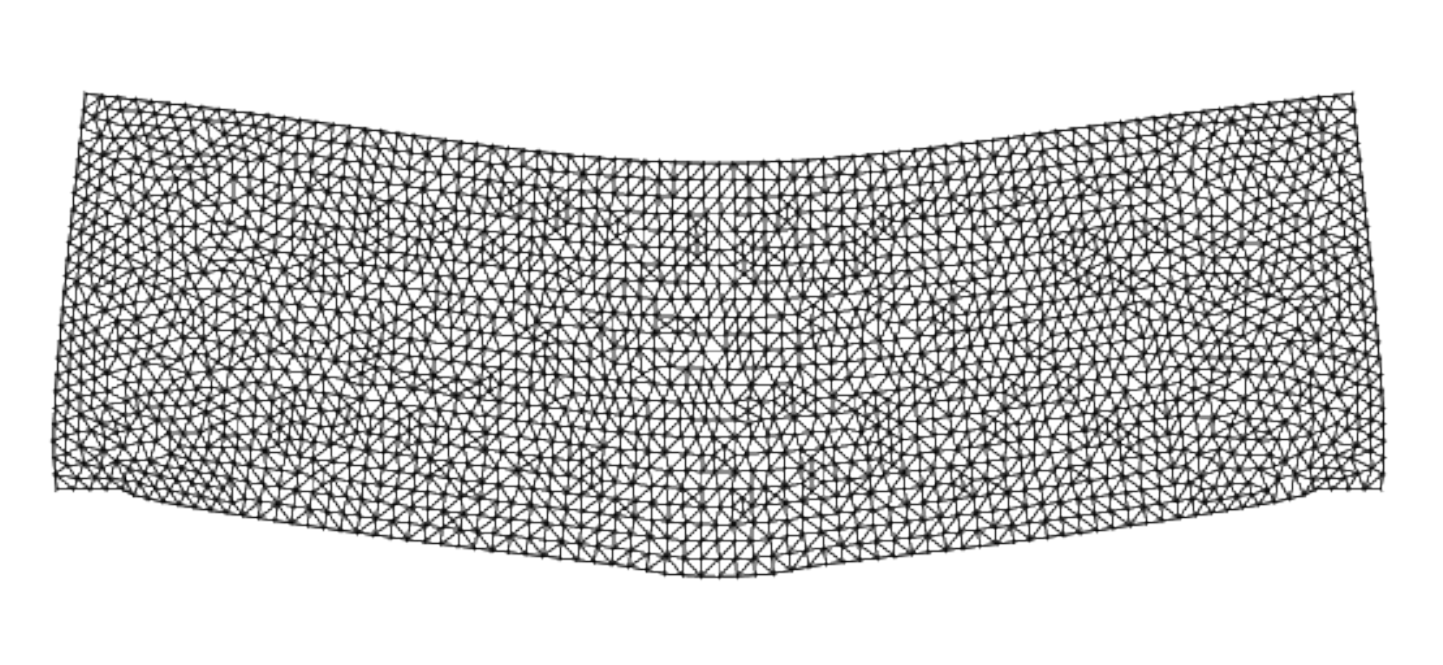}\quad     
\includegraphics[width=6cm]{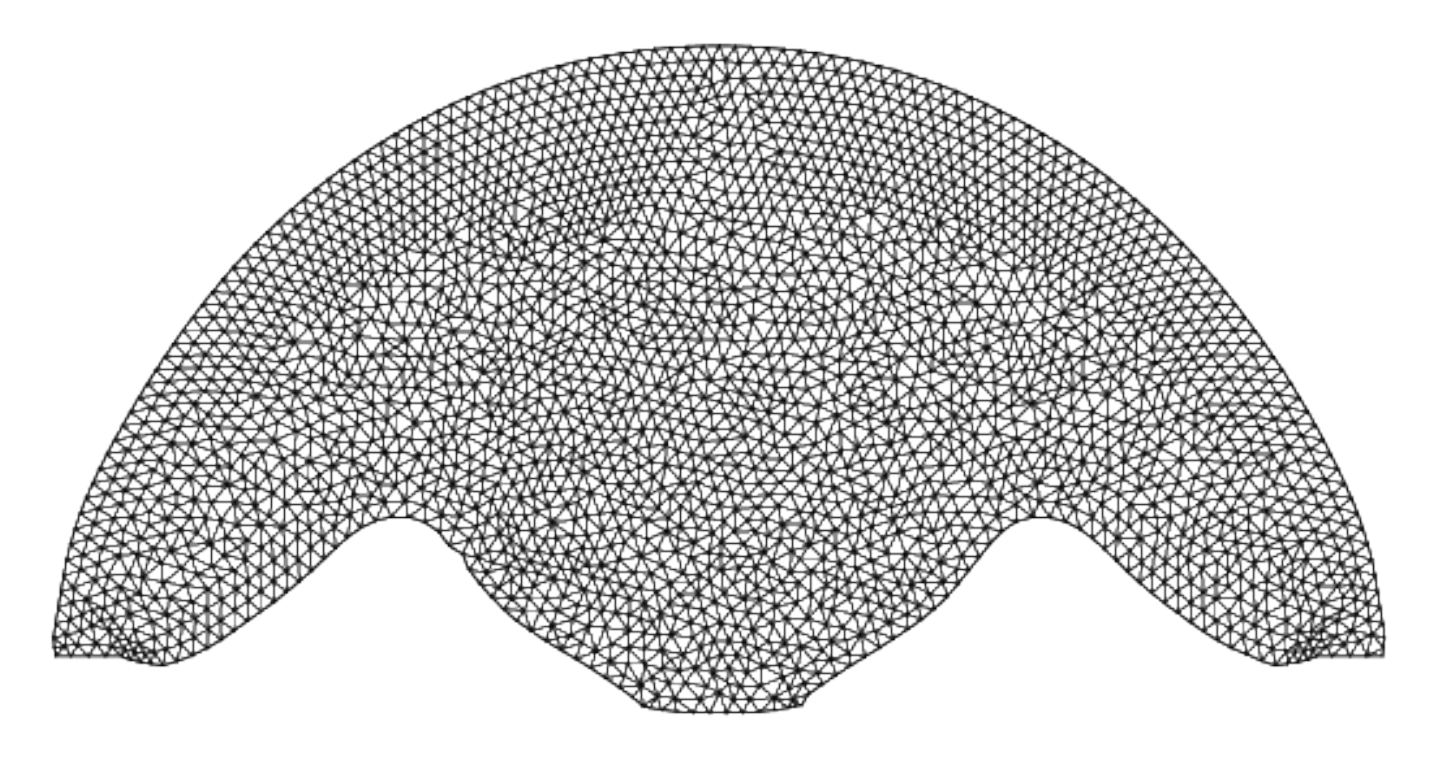}
\end{center}
  \caption{Bridge. The meshes after the elastic deformations for the initial domain
    $\Omega_0=]-1,1[ \times ]0,0.6[$ and
    for the optimal domain $\Omega_{100}$ presented in Figure \ref{fig:pont2_sigma_half_D}.
  The displacements were reduced by a factor $0.1$.
The cost (\ref{compliance1}) decreases from $0.378632$ (left image) 
to $0.297857$ (right image).}
\label{fig:pont2_mesh}
\end{figure}

Furthermore, we have used for the initialization $\Omega_0=]-1,1[ \times ]0,0.6[$ obtained for
$g_0(x_1,x_2)=0.1(0.6-x_2)$.
The initial value of the objective function
is $J(g_0)=0.353644$ and the value after 100 iterations  (\textbf{Step 5}) is
$J(g_{100})=0.296596$. The optimal domain is presented in Figure \ref{fig:pont2_sigma_half_D}.
We have also solved the original elasticity problem (\ref{elast1})-(\ref{elast4}), for the
initial and final domains.
Using \texttt{FreeFem++}, it is possible to build a mesh which boundary is the zero level set
of a function $g$. From Proposition \ref{propNEW:1}, the solution computed in a mesh of
$\Omega_g$ is close to the solution
of (\ref{stateeps}) computed in a fixed mesh of $D$.
The deformations are presented in Figure \ref{fig:pont2_mesh}.
The values of the cost (\ref{compliance1}) are $0.378632$ for the initial domain
and $0.297857$ for the final domain.

The above {\color{black}final domains, that are not necessarily globally optimal, }
  differ from the solutions obtained using homogenization or level set methods,
\cite{Allaire2004}, \cite{Allaire2014}. This is no contradiction since the considered
optimization problems are highly non convex and the solutions depend on the initial
iteration, or on the chosen parameters, etc. The important characteristic is that the
method discussed in this paper ensures a consistent decrease of the cost together
with topological and boundary variations for the obtained domains.


Finally, we have also tested the dependence on $\epsilon$ of the state solution of (\ref{stateeps})
and Prop. \ref{propNEW:1} was confirmed.
The errors in the norms $L^2$ and $H^1$ are reported in Table \ref{tab:pont_init}
for the initial domain $\Omega_0=]-1,1[ \times ]0,0.6[$ and in Table \ref{tab:pont_fin}
for the final domain $\Omega_{100}$. The reference displacement $\mathbf{y}^*$ is the
solution of (\ref{stateeps}) with $\widehat{H}$ in place of $H^\epsilon$, where $\widehat{H}(r)$
takes 1 for $r\geq 0$ and it takes $10^{-9}$ for $r<0$.
\begin{table}[ht]
\begin{center}
\begin{tabular}{|l|c|c|c|}
\hline
$\epsilon $ &
$J$ &
$\left\| \mathbf{y}^\epsilon(g_0)-\mathbf{y}^*\right\|_{L^2(\Omega_0)}$ &
$\left\| \mathbf{y}^\epsilon(g_0)-\mathbf{y}^*\right\|_{H^1(\Omega_0)}$\\  \hline
0.01   & 0.353644 & 0.097841 & 0.302369
\\  \hline
0.005  & 0.369480 & 0.035030 & 0.112965
\\  \hline
0.001  & 0.378150 & 0.002096 & 0.027376
\\  \hline
0.0005 & 0.378506 & 0.000799 & 0.026536
\\  \hline
\end{tabular}
\end{center}
\caption{The dependence of the cost given by (\ref{compliance1}) and of the displacement
  $\mathbf{y}^\epsilon(g_0)$ solution of (\ref{stateeps}) in the initial domain
  $\Omega_0=]-1,1[ \times ]0,0.6[$. The cost for $\mathbf{y}^*$ is $0.378727$.}
\label{tab:pont_init}
\end{table}

\clearpage
\begin{table}[ht]
\begin{center}
\begin{tabular}{|l|c|c|c|}
\hline
$\epsilon $ &
$J$ &
$\left\| \mathbf{y}^\epsilon(g_{100})-\mathbf{y}^*\right\|_{L^2(\Omega_{100})}$ &
$\left\| \mathbf{y}^\epsilon(g_{100})-\mathbf{y}^*\right\|_{H^1(\Omega_{100})}$\\  \hline
0.01   & 0.296596 & 0.018279 & 0.105998
\\  \hline
0.005  & 0.297032 & 0.005071 & 0.084750
\\  \hline
0.002  & 0.297813 & 0.002263 & 0.076559
\\  \hline
0.001  & 0.298063 & 0.001925 & 0.073167
\\  \hline
\end{tabular}
\end{center}
\caption{The dependence of the cost given by (\ref{compliance1}) and of the displacement
  $\mathbf{y}^\epsilon(g_{100})$ solution of (\ref{stateeps}) in the final domain
  $\Omega_{100}$. The cost for $\mathbf{y}^*$ is $0.298536$.}
\label{tab:pont_fin}
\end{table}

The differences between $\mathbf{y}^*$ and the original elasticity problem (\ref{elast1})-(\ref{elast4})
are $0.000362$ and $0.057692$ in the norms $L^2$ and $H^1$ respectively for the initial domain $\Omega_0$.
The cost difference is $|0.378727-0.378632|=0.000095$.
For the final domain  $\Omega_{100}$, the differences are $0.002925$ and $0.110506$ in the norms
$L^2$ and $H^1$ respectively
and the cost difference is $|0.298536-0.297857|=0.000679$. We notice that the descent property
remains valid in the original shape optimization problem as well, due to the good approximation
provided by our fixed domain approach.

\end{document}